\documentclass[12pt]{amsart}

\usepackage{amsmath}
\usepackage{amssymb}
\usepackage{amsfonts}

\usepackage[foot]{amsaddr}

\usepackage{subfigure}
\usepackage{epsfig}

\usepackage{graphicx}

\newtheorem{thm}{Theorem}[section]
\newtheorem{dfn}[thm]{Definition}
\newtheorem{cor}[thm]{Corollary}
\newtheorem{prop}[thm]{Proposition}
\newtheorem{lem}[thm]{Lemma}
\theoremstyle{remark} 
\newtheorem{rmk}[thm]{Remark}

\def\cqfd{\mbox{}\nolinebreak\hfill$\Box$\medbreak\par}
\newenvironment{pf}{\noindent\textbf{Proof:}}{\cqfd}

\newcommand{\B}{\mathcal{B}}
\newcommand{\hB}{\widehat{\mathcal{B}}}
\newcommand{\OB}{\overline{\mathcal{B}}}

\newcommand{\C}{\mathbb{C}}

\newcommand{\R}{\mathbb{R}}
\newcommand{\Q}{\mathbb{Q}}
\def\rh#1{{{\bf H}^{#1}_{\R}}}
\def\ch#1{{{\bf H}^{#1}_{\C}}}
\def\cp#1{{{\bf P}^{#1}_{\C}}}
\def\chb#1{{\overline{\bf H}^{#1}_{\C}}}
\def\pu21{\mathbf{PU}(2,1)}
\def\un21{\mathbf{U}(2,1)}

\usepackage[OT2,OT1]{fontenc}
\newcommand\cyr{%
\renewcommand\rmdefault{wncyr}%
\renewcommand\sfdefault{wncyss}%
\renewcommand\encodingdefault{OT2}%
\normalfont
\selectfont}
\DeclareTextFontCommand{\textcyr}{\cyr}

\title{Complex hyperbolic geometry of the figure eight knot}

\begin{document}

\author{Martin Deraux}
\author{Elisha Falbel}

\date{\today}

\begin{abstract}
  We show that the figure eight knot complement admits a uniformizable
  spherical CR structure, i.e. it occurs as the manifold at infinity
  of a complex hyperbolic orbifold. The uniformization is unique
  provided we require the peripheral subgroups to have unipotent
  holonomy.
\end{abstract}

\maketitle

\section{Introduction} \label{sec:intro}

The general framework of this paper is the study of the interplay
between topological properties of $3$-manifolds and the existence of
geometric structures. The model result along these lines is of course
Thurston's geometrization conjecture, recently proved by Perelman,
that contains a topological characterization of manifolds that admit a
geometry modeled on real hyperbolic space $\rh 3$. Beyond an existence
result (under the appropriate topological assumptions), the hyperbolic
structures can in fact be constructed fairly explicitly, as one can
easily gather by reading Thurston's notes~\cite{thurstonNotes}, where
a couple of explicit examples are worked out. 

The idea is to triangulate the manifold, and to try and realize each
tetrahedron geometrically in $\rh 3$. The gluing pattern of the
tetrahedra imposes compatibility conditions on the parameters of the
tetrahedra, and it turns out that solving these compatibility
equations is very often equivalent to finding the hyperbolic
structure. The piece of software called SnapPea, originally developed
by Jeff Weeks (and under constant development to this day), provides
an extremely efficient way to construct explicit hyperbolic structures on
3-manifolds.

In this paper, we are interested in using the $3$-sphere $S^3$ as the
model geometry, with the natural structure coming from describing it
as the boundary of the unit ball $\mathbb{B}^2\subset\C^2$. Any real
hypersurface in $\C^2$ inherits what is called a CR structure (the
largest subbundle in the tangent bundle that is invariant under the
complex structure), and such a structure is called spherical when it
is locally equivalent to the CR structure of $S^3$. Local equivalence
to $S^3$ in the sense of CR structures translates into the existence
of an atlas of charts with values in $S^3$, and with transition maps
given by restrictions of biholomorphisms of $\mathbb{B}^2$, i.e. elements of
$\pu21$, see~\cite{burnsshnider}.

In other words, a spherical CR structure is a $(G,X)$-structure with
$G=\pu21$, $X=S^3$. The central motivating question is to give a
characterization of $3$-manifolds that admit a spherical CR
structure; the only negative result in that direction is given by
Goldman~\cite{goldmanNilpotent}, who classifies $T^2$-bundles over
$S^1$ that admit spherical CR structures (only those with Nil geometry
admit spherical CR structures).

An important class of spherical CR structures is the class of
\emph{uniformizable} spherical CR structures. These are obtained from
discrete subgroups $\Gamma\subset \pu21$ by taking the quotient of the
domain of discontinuity $\Omega$ by the action of $\Gamma$ (we assume
that $\Omega$ is non-empty, and that $\Gamma$ has no fixed point on
$\Omega$, so that the quotient is indeed a manifold). The structure induced
from the standard CR structure on $S^3$ on the quotient
$M=\Gamma\setminus\Omega$ is then called a \emph{uniformizable
  spherical CR structure} on $M$.

When a manifold $M$ can be written as above for some group $\Gamma$,
we will also simply say that $M$ admits a spherical CR
uniformization. Our terminology differs slightly from the recent
literature on the subject, where uniformizable structures are
sometimes referred to as \emph{complete} structures
(see~\cite{richBook} for instance).

Of course one wonders which manifolds admit spherical CR
uniformizations, and how restrictive it is to require the existence of
a spherical CR uniformization as opposed to a general spherical CR
structure. For instance, when $\Gamma$ is a finite group acting
without fixed points on $S^3$, $\Omega=S^3$ and $\Gamma\setminus S^3$
gives the simplest class of examples (including lens spaces).

The class of circle bundles over surfaces has been widely explored,
and many such bundles are known to admit uniformizable spherical CR
structures, see the introduction of~\cite{richBook} and the references
given there. It is also known that well-chosen deformations of
triangle groups produce spherical CR structures on more complicated
$3$-manifolds, including real hyperbolic ones. Indeed, Schwartz showed
in~\cite{richWhitehead} that the Whitehead link complement admits a
uniformizable spherical CR structure, and in~\cite{schwartz4447} he
found an example of a closed hyperbolic manifold that arises as the
boundary of a complex hyperbolic surface. Once again, we refer the
reader to the~\cite{richBook} for a detailed overview of the history
of this problem.
 
All these examples are obtained by analyzing special classes of
discrete groups, and checking the topological type of their manifold
at infinity. In the opposite direction, given a $3$-manifold $M$, one
would like a method to construct (and possibly classify) all
structures on $M$, in the spirit of the constructive version of
hyperbolization alluded to earlier in this introduction.

A step in that direction was proposed by the second author
in~\cite{falbelFigure8}, based on triangulations and adapting the
compatibility equations to the spherical CR setting. Here, a basic
difficulty is that there is no canonical way to associate a
tetrahedron to a given quadruple of points in $S^3$. Even the
$1$-skeleton is elusive, since arcs of $\C$-circles (or $\R$-circles)
between two points are not unique (see section~\ref{sec:complex} for
definitions).

A natural way over this difficulty is to formulate compatibility
conditions that translate the possibility of geometric realization in
$S^3$ only on the level of the vertices of the tetrahedra. Indeed,
ordered generic quadruples of points are parametrized up to isometry
by appropriate cross ratios, and one can easily write down the
corresponding compatibility conditions explicitly~\cite{falbelFigure8}.

Given a solution of these compatibility equations, one always gets a
representation $\rho:\pi_1(M)\rightarrow \pu21$, but it is not clear
whether or not the quadruples of points can be extended to actual
tetrahedra in a $\rho$-equivariant way (in other words, it is not
clear whether or not $\rho$ is the holonomy of an actual structure).

There are many solutions to the compatibility equations,
so we will impose a restriction on the representation $\rho$, namely
that $\rho(\pi_1(T))$ be unipotent for each torus boundary component
$T$ of $M$. This is a very stringent condition, but it is natural
since it holds for complete hyperbolic metrics of finite volume.

For the remainder of the paper, we will concentrate on a specific
$3$-manifold, namely the figure eight knot complement, and give
encouraging signs for the philosophy outlined in the preceding
paragraphs. Indeed, for that specific example, we will check that the
solutions to the compatibility equations give a spherical CR
uniformization of the figure eight knot, which is unique provided we
require the boundary holonomy to consist only of unipotent isometries
(in fact we get one structure for each orientation on $M$, see
section~\ref{sec:mcg}).

We work with the figure eight knot complement partly because it played
an important motivational role in the eighties for the development of
real hyperbolic geometry. It is well known that this non-compact
manifold $M$ admits a unique complete hyperbolic metric, with one
torus end (which one may think of as a tubular neighborhood of the
figure eight knot). This is originally due to Riley,
see~\cite{rileyfigure8}.

\begin{figure}[htbp]
  \centering
  \epsfig{figure=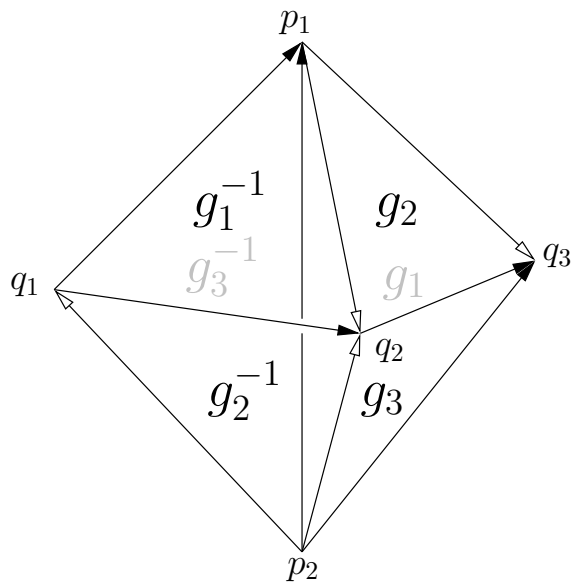, width=0.6\textwidth}
  \caption{The figure eight knot complement can be obtained by gluing
    two tetrahedra (a face on the left and a face on the right are
    identified if the corresponding pattern of arrows agree), and
    removing the vertices.}\label{fig:tetrahedra}
\end{figure}

It is also well known that $M$ can be triangulated with just two
tetrahedra (this triangulation is far from simplicial, but this is
irrelevant in the present context). The picture in
Figure~\ref{fig:tetrahedra} can be found for instance in the first few
pages of Thurston's notes~\cite{thurstonNotes}.
The above decomposition can be realized geometrically in $\rh 3$ (and
the corresponding geometric tetrahedra are regular tetrahedra, so the
volume of this metric is $6\textcyr{L}(\pi/3)\approx 2.029$).

For the specific triangulation of the figure eight knot complement
depicted in Figure~\ref{fig:tetrahedra}, all the solutions of the
compatibility equations were given in~\cite{falbelFigure8}, without
detailed justification of the fact that the list of solutions is
exhaustive. The explanation of exhaustivity now appears in various
places in the literature (see~\cite{bfg} and~\cite{ggz}, and
also~\cite{censusfkr} for more general 3-manifolds). It
turns out there are only three solutions to the compatibility
equations (up to complex conjugation of the cross ratios parametrizing
the tetrahedra), yielding three representations $\rho_1$, $\rho_2$ and
$\rho_3:\pi_1(M)\rightarrow \pu21$ (in fact six representations, if we
include their complex conjugates). Throughout the paper, we will
denote by $\Gamma_k$ the image of $\rho_k$.

It was shown in~\cite{falbelFigure8} that $\rho_1$ is the holonomy of
a \emph{branched} spherical CR structure (the corresponding developing
map is a local diffeomorphism away from a curve), and that the limit
set of $\Gamma_1$ is equal to $\partial_\infty \ch 2$, hence the
quotient $\Gamma_1\setminus \ch 2$ has empty manifold at infinity. In
particular, no spherical CR structure with holonomy $\rho_1$ can ever
be uniformizable. In~\cite{falbelwang}, a branched structure with
holonomy $\rho_2$ is constructed, which is again not a uniformization.

The main goal of this paper is to show that $\rho_2$ and $\rho_3$ are
holonomy representations of unbranched uniformizable spherical CR
structures on the figure eight knot complement. These two
representations are not conjugate in $\pu21$, but it turns out that
the images $\Gamma_2$ and $\Gamma_3$ are in fact conjugate.

The precise relationship between the two structures corresponding to
$\rho_2$ and $\rho_3$ will be explained by the existence of an
orientation-reversing diffeomorphism of the figure eight knot
complement (which follows from the fact that this knot is
amphichiral). Indeed, given a diffeomorphism $\alpha:M\rightarrow
\Gamma_2\setminus \Omega_2$, and an an orientation-reversing
diffeomorphism $\varphi:M\rightarrow M$, $\alpha\circ\varphi$ defines a
spherical CR structure on $M$ with the opposite orientation.  We will
see that $\rho_2$ and $\rho_3$ are obtained from each other by this
orientation switch (see section~\ref{sec:thirdRep}). For that reason,
we will work only with $\rho_2$ for most of the paper.

We denote by $\Gamma$ the group $\Gamma_2$. Our main result is the
following.
\begin{thm}\label{thm:main}
 The domain of discontinuity $\Omega$ of $\Gamma$ is non empty. The
 action of $\Gamma$ has no fixed points in $\Omega$, and the quotient
 $\Gamma\setminus\Omega$ is homeomorphic to the figure eight knot
 complement.
\end{thm}
In other words, the figure eight knot admits a
spherical CR uniformization, with uniformization given by
$\Gamma$. The uniformization is not quite unique, but we will show
that it is unique provided we require the boundary holonomy to be
unipotent (see Proposition~\ref{prop:unique}).

The fact that the ideal boundary of $\Gamma\setminus \ch 2$ is
indeed a manifold, and not just an orbifold, follows from the fact
that every elliptic element in $\Gamma$ has an isolated fixed point in
$\ch 2$ (we will be able to list all conjugacy classes of elliptic
elements, by using the cycles of the fundamental domain, see
Proposition~\ref{prop:nofixedpoint}).

The result of Theorem~\ref{thm:main} is stated in terms of the domain
of discontinuity which is contained in $\partial_\infty \ch 2$, so one
may expect the arguments to use properties of $S^3\subset \C^2$ or
Heisenberg geometry (see section~\ref{sec:complex}). In fact the bulk of the
proof is about the relevant complex hyperbolic orbifold
$\Gamma\setminus \ch 2$, and for most of the paper, we will use
geometric properties of $\ch 2$.

The basis of our study of the manifold at infinity will be the
Dirichlet domain for $\Gamma$ centered at a strategic point, namely
the isolated fixed point of $G_2=\rho_2(g_2)$ (see
section~\ref{sec:threereps} for notation). This domain is not a
fundamental domain for the action of $\Gamma$ (the center is
stabilized by a cyclic group of order $4$), but it is convenient
because it has very few faces (in fact all its faces are isometric to
each other). In particular, we get an explicit presentation for
$\Gamma$, given by
\begin{equation}\label{eq:presentation}
\langle\ G_1, G_2\ |\ G_2^4,\ (G_1G_2)^3,\ (G_2G_1G_2)^3\ \rangle
\end{equation}
Note that $\rho_2$ is of course not a faithful representation of the
figure eight knot group. In fact from the above presentation, it is
easy to determine normal generators for the kernel of $\rho_2$, see
Proposition~\ref{prop:kernel}.\\

\noindent {\bf Acknowledgements:} This work was supported in part by
the ANR through the project Structure G\'eom\'etriques et
Triangulations. The authors wish to thank Christine Lescop, Antonin
Guilloux, Julien March\'e, Jean-Baptiste Meilhan, Jieyan Wang, Pierre
Will and Maxime Wolff for useful discussions related to the
paper. They are also very grateful to the referees for many useful
remarks that helped improve the exposition.

\section{Basics of complex hyperbolic geometry}

\subsection{Complex hyperbolic geometry} \label{sec:complex}

In this section we briefly review basic facts and notation about the
complex hyperbolic plane. For more information,
see~\cite{goldman}.

We denote by $\C^{2,1}$ the three-dimensional complex vector space
$\C^3$ equipped with the Hermitian form
$$
\langle Z, W \rangle=Z_1\overline{W}_3+Z_2\overline{W}_2+Z_3\overline{W}_1.
$$

The subgroup of $GL(3,\C)$ preserving the Hermitian form $\langle
\cdot,\cdot \rangle$ is denoted by $\un21$, and its action preserves
each of the following three sets:
$$
  V_{+}=  \{ Z\in \mathbb{C}^{2,1} : \langle Z, Z \rangle >0 \},
$$
$$
   V_{0}=  \{ Z\in \mathbb{C}^{2,1}- \{0\} : \langle Z, Z \rangle =0 \},
$$
$$
   V_{-}=  \{ Z\in \mathbb{C}^{2,1} : \langle Z, Z \rangle < 0 \}.
$$ 

Let $P : \mathbb{C}^{2,1}- \{0\} \rightarrow \cp 2 $ be the canonical
projection onto complex projective space, and let $\pu21$ denote the
quotient of $\un21$ by scalar matrices, which acts
effectively on $\cp 2$. Note that the action of $\pu21$ is transitive
on $P(V_{\pm})$ and on $P(V_0)$. Up to scalar multiples, there is a
unique Riemannian metric on $P(V_-)$ invariant under the action of
$\pu21$, which turns it into a Hermitian symmetric space often denoted
by $\ch 2$, and called the complex hyperbolic plane. In the present
paper, we will not need a specific normalization of the metric. We
mention for completeness that any invariant metric is K\"ahler, with
holomorphic sectional curvature a negative constant (the real
sectional curvatures are $1/4$-pinched).

The full isometry group of $\ch 2$ is given by 
$$
\widehat{\pu21}=\langle \pu21, {\mathbf\iota} \rangle,
$$ where ${\mathbf\iota}$ is given by complex conjugation $Z\mapsto
\overline{Z}$ on the level of homogeneous coordinates.

Still denoting $(Z_1,Z_2,Z_3)$ the coordinates of $\C^3$, one easily
checks that $V_-$ can contain no vector with $Z_3=0$, hence we can
describe its image in $\cp 2$ in terms of non-homogeneous coordinates
$w_1=Z_1/Z_3$, $w_2=Z_2/Z_3$, where $P(V_-)$ corresponds to the Siegel
half space
$$
|w_1|^2+2\ {\rm Re}\ w_2<0.
$$ 
The ideal boundary of complex hyperbolic space is defined as
$\partial_\infty \ch 2=P(V_{0})$. It is described almost entirely in
the affine chart $Z_3\neq 0$ used to define the Siegel half space,
only $(1,0,0)$ is sent off to infinity. We denote by $p_\infty$ the
corresponding point in $\partial_\infty \ch2$.

The unipotent stabilizer of $(1,0,0)$ acts simply transitively on
$\partial_\infty \ch 2\setminus\{p_\infty\}$, which allows us to
identify $\partial_\infty\ch 2$ with the one-point compactification of
the Heisenberg group $\mathfrak{N}$.

Here recall that $\mathfrak{N}$ is defined as $\mathbb{C}\times
\mathbb{R}$ equipped with the following group law
$$
(z,t)\cdot (z',t')=(z+z',t+t'+2\Im(z\overline{z}')).
$$ 
Any point $p=(z,t)\in \mathfrak{N}$ has the following lift to
$\mathbb{C}^{2,1}$:
$$
\tilde{p}=\left[\begin{array}{c}
          (-|z|^2+it)/2 \\
          z \\
          1
        \end{array}\right]
$$
while $p_\infty$ lifts to $(1,0,0)$.

It is a standard fact that the above form can be diagonalized, say by
using the change of homogeneous coordinates given by $U_2=Z_2$,
$U_1=(Z_1+Z_3)/\sqrt{2}$, $U_3=(Z_1-Z_3)/\sqrt{2}$. With these
coordinates, the Hermitian form reads
$$
\langle U,V\rangle = U_1\overline{V}_1 + U_2\overline{V}_2 - U_3\overline{V}_3, 
$$ 
and in the affine chart $U_3\neq 0$, with coordinates $u_1=U_1/U_3$,
$u_2=U_2/U_3$, $\ch 2$ corresponds to the unit ball $\mathbb{B}^2\subset
\C^2$, given by
$$
|u_1|^2+|u_2|^2<1.
$$ 
In this model the ideal boundary is simply given by the unit sphere
$S^3\subset \C^2$. This gives $\partial_\infty \ch 2$ a natural
CR-structure (see the introduction and the references given there).

We will use the classification of isometries of negatively curved
spaces into elliptic, parabolic and loxodromic elements, as well as a
slight algebraic refinement; an elliptic isometry is called {\bf
  regular elliptic} if its matrix representatives have distinct
eigenvalues.

Non-regular elliptic elements in $\pu21$ fix a projective line in
$\cp 2$, hence they come into two classes, depending on the position
of that line with respect to $\ch 2$. If the projective line
intersects $\ch 2$, the corresponding isometry is called a {\bf
  complex reflection in a line}; if it does not intersect
$\partial_\infty \ch2$, then the isometry is called a {\bf complex
  reflection in a point}. Complex reflections in points do not have any
fixed points in the ideal boundary.

The only parabolic elements we will use in this paper will be {\bf
  unipotent} (i.e. some matrix representative in $\un21$ has $1$ as
its only eigenvalue).

Finally, we mention the classification of totally geodesic
submanifolds in $\ch 2$. There are two kinds of totally geodesic
submanifolds of real dimension two, complex geodesics (which can be
thought of copies of $\ch 1$), and totally real totally geodesic
planes (copies of $\rh 2$).

In terms of the ball model, complex lines correspond to intersections
with $\mathbb{B}^2$ of affine lines in $\C^2$. In terms of projective
geometry, they are parametrized by their so-called polar vector, which
is the orthogonal complement of the corresponding plane in $\C^3$ with
respect to the Hermitian form $\langle \cdot,\cdot\rangle$.

The trace on $\partial_\infty \ch2$ of a complex geodesic (resp. of a
totally real totally geodesic plane) is called a $\C$-circle (resp. an
$\R$-circle). 

For completeness, we mention that there exists a unique complex line
through any pair of distinct points $p,q\in\partial_\infty \ch 2$. The
corresponding $\C$-circle is split into two arcs, but there is in
general no preferred choice of an arc of $\C$-circle between $p$ and
$q$. Given $p,q$ as above, there are infinitely many $\R$-circles
containing them. The union of all these $\R$-circles is called a
spinal sphere (see section~\ref{sec:bisectors} for more on this).

\subsection{Generalities on Dirichlet domains} \label{sec:dirichlet}

Recall that the Dirichlet domain for $\Gamma\subset \pu21$ centered
at $p_0\in\ch 2$ is defined as
$$ 
E_\Gamma = \left\{ z\in \ch 2 : d(z,p_0)\leqslant d(z,\gamma p_0)
\textrm{ for all }\gamma\in \Gamma\right\}.
$$ 
Although this infinite set of inequalities is in general quite hard
to handle, in many situations there is a finite set of inequalities
that suffice to describe the same polytope (in other words, the
polytope has finitely many faces).

Given a (finite) subset $S\subset \Gamma$, we denote by 
$$ 
E_S = \left\{ z\in \ch 2 : d(z,p_0)\leqslant d(z,\gamma p_0) \textrm{ for
  all }\gamma\in S\right\},
$$ 
and search for a minimal set $S$ such that $E_\Gamma=E_S$. In
particular, we shall always assume that
\begin{itemize}
  \item $s p_0\neq p_0$ for every $s\in S$ and
  \item $s_1p_0\neq s_2 p_0$ for every $s_1\neq s_2\in S$.
\end{itemize} 
Indeed, $sp_0=p_0$ would give a vacuous inequality, and $s_1 p_0=s_2
p_0$ would give a repeated face.

Given a finite set $S$ as above and an element $\gamma\in S$, we refer
to the set of points equidistant from $p_0$ and $\gamma p_0$ as {\bf
  the bisector associated to $\gamma$}, i.e.
$$ 
\B(p_0,\gamma p_0)=\left\{z\in \ch 2: d(z,p_0)=d(z,\gamma
p_0)\right\}.
$$ 
We will say that $\gamma$ {\bf defines a face} of $E_S$ when
$\B(p_0,\gamma p_0)\cap E_S$ has non empty interior in $\B(p_0,\gamma
p_0)$. In that case, we refer to $\B(p_0,\gamma p_0)\cap E_S$ as {\bf
  the face of $E_S$ associated to $\gamma$}.

We will index the bisectors bounding $E_S$ by integers $k$, and write
$\B_k$ for the $k$-th bounding bisector. We will then often write
$b_k$ for the corresponding face, i.e. $b_k=\B_k\cap E_S$ (this
notation only makes sense provided the set $S$ is clear from the
context, which will be the case later in the paper).

The precise determination of all the faces of $E_S$, or equivalently
the determination of a minimal set $S$ with $E_S=E_\Gamma$ is quite
difficult in general. 

The main tool for proving that $E_\Gamma=E_S$ is the Poincar\'e
polyhedron theorem, which gives sufficient conditions for $E_S$ to be a
fundamental domain for the group generated by $S$. The assumptions are
roughly as follows:
\begin{enumerate}
  \item $S$ is symmetric (i.e. $\gamma^{-1}\in S$ whenever $\gamma\in
    S$) and the faces of $E_S$ associated to $\gamma$ and $\gamma^{-1}$ are
    isometric.
  \item The images of $E_S$ under elements of $\Gamma$ give a local
    tiling of $\ch 2$.
\end{enumerate}
The conclusion of the Poincar\'e polyhedron theorem is then that the
images of $E_S$ under the group generated by $S$ give a global tiling
of $\ch 2$ (from this one can deduce a presentation for the group
$\langle S\rangle$ generated by $S$).

The requirement that opposite faces be isometric justifies calling the
elements of $S$ ``side pairings''. We shall use a version of the
Poincar\'e polyhedron theorem for \emph{coset decompositions} rather
than for groups, because we want to allow some elements of $\Gamma$ to
fix the center $p_0$ of the Dirichlet domain.

The result we have in mind is stated for the simpler case of $\ch 1$
in~\cite{beardon}, section~9.6. We assume $E_S$ is stabilized by a
certain (finite) subgroup $H\subset \Gamma$, and the goal is to show
that $E_S$ is a fundamental domain modulo the action of $H$, i.e. if
$\gamma_1 E_S\cap \gamma_2 E_S$ has non empty interior, then
$\gamma_1=\gamma_2 h$ for some $h\in H$.

The corresponding statement for $\ch 2$ appears
in~\cite{mostowPacific}, with a light treatment of the assumptions
that guarantee completeness, so we list the hypotheses roughly as they
appear in~\cite{maskit} (see also~\cite{JRP-book} for a proof in the
context of complex hyperbolic space). The local tiling condition will
consist of two checks, one for ridges (faces of codimension two in
$E_S$), and one for boundary vertices.  A ridge $e$ is given by the
intersection of two faces of $E_S$, i.e. two elements $s,t\in S$. We
will call the intersection of $E_S$ with a small tubular neighborhood
of $e$ the wedge of $E_S$ near $e$.

\begin{itemize}
\item Given a ridge $e$ defined as the intersection of two faces
  corresponding to $s,t\in S$, we consider all the other ridges of
  $E_S$ that are images of $e$ under successive side pairings or
  elements of $H$, and check that the corresponding wedges tile a
  neighborhood of that ridge.

\item Given a boundary vertex $p$, which is given by (at least) three
  elements $s,t,u\in S$, we need to consider the orbit of $p$ in $E_S$
  using successive side pairings or elements of $H$, check that the
  corresponding images of $E_S$ tile a neighborhood of that vertex,
  and that the corresponding cycle transformations are all given by
  \emph{parabolic} isometries.
\end{itemize}
 
The conclusion of the Poincar\'e theorem is that if $\gamma_1 E_S\cap
\gamma_2 E_S$ has non-empty interior, then $\gamma_1$ and $\gamma_2$
differ by right multiplication by an element of $H$. From this, one
easily deduces a presentation for $\Gamma$, with generators given by
$S\cup H$ ($H$ can of course be replaced by any generating set for
$H$), and relations given by ridge cycles (together with the relations
in a presentation of $H$).

\subsection{Bisector intersections}\label{sec:bisectors}

In this section, we review some properties of bisectors and bisector
intersections (see~\cite{goldman} or~\cite{deraux4445} for more
information on this).

Let $p_0,p_1\in \ch 2$ be distinct points given in homogeneous
coordinates by vectors $\tilde{p}_0$, $\tilde{p}_1$, chosen so that
$\langle \tilde{p}_0,\tilde{p}_0 \rangle=\langle
\tilde{p}_1,\tilde{p}_1 \rangle$. By definition, the {\bf bisector}
$\B=\B(p_0,p_1)$ is the locus of points equidistant of $p_0,p_1$. It
is given in homogeneous coordinates $\mathbf{z}=(z_0,z_1,z_2)$ by the
negative vectors $\mathbf{z}$ that satisfy the equation
\begin{equation}\label{eq:bisector}
|\langle \mathbf{z},\tilde{p}_0\rangle|=|\langle
\mathbf{z},\tilde{p}_1\rangle|.
\end{equation}
When $\mathbf{z}$ is not assumed to be negative, the same equation
defines an {\bf extor} in projective space. Note that $\mathbf{z}$ is
a solution to this equation if and only if it is orthogonal (with
respect to the indefinite Hermitian inner product) to some vector of
the form $\tilde{p}_0-\alpha\tilde{p}_1$, with $|\alpha|=1$.

Finally, we mention that the image in projective space of the set of
null vectors $\mathbf{z}$, i.e. such that $\langle
\mathbf{z},\mathbf{z}\rangle=0$, and that satisfy
equation~(\ref{eq:bisector}) is a topological sphere, which we will
call either the {\bf boundary at infinity} corresponding to the
bisector, or its {\bf spinal sphere}.

Restricting to vectors $\tilde{p}_0-\alpha\tilde{p}_1$ which have
positive square norm, we get a foliation of $\B(p_0,p_1)$ by complex
lines given by the set of negative lines in
$(\tilde{p}_0-\alpha\tilde{p}_1)^\perp$ for fixed value of
$\alpha$. These complex lines are called the {\bf complex slices} of
the bisector. Negative vectors of the from
$(\tilde{p}_0-\alpha\tilde{p}_1)$ (still with $|\alpha|$=1)
parametrize a real geodesic, which is called the {\bf real spine} of
$\B$. The complex geodesic that it spans is called the {\bf complex
  spine} of $\B$.  There is a natural extension of the real spine to
projective space, given by the (not necessarily negative) vectors of the
form $\tilde{p}_0-\alpha\tilde{p}_1$, we call this the {\bf extended
  real spine} (the complex projective line that contains it is called
the {\bf extended complex spine}).

Geometrically, each complex slice of $\B$ is the preimage of a given
point of the real spine under orthogonal projection onto the complex
spine, and in particular, the bisector is uniquely determined by its
real spine. 

Given two distinct bisectors $\B_1$ and $\B_2$, their intersection is
to a great extent controlled by the respective positions of their
complex spines $\Sigma_1$ and $\Sigma_2$. In particular, if $\Sigma_1$
and $\Sigma_2$ intersect outside of their respective real spines, the
bisectors are called {\bf coequidistant}.

This special case of bisector intersections is important in the
context of Dirichlet domains, since by construction all the faces of a
Dirichlet domain are equidistant from one given point (namely its
center). We recall the following, which is an important tool for
studying the combinatorics of polyhedra bounded by bisectors (and also
in order to apply the Poincar\'e polyhedron theorem, see
section~\ref{sec:poincare}).
\begin{thm}\label{thm:giraud}
  Let $\B_1$ and $\B_2$ be coequidistant bisectors. Then their
  intersection is a smooth disk, which is contained in precisely three
  bisectors.
\end{thm}
This theorem is due to Giraud (for a detailed proof see sections~8.3.5
and~9.2.6 of~\cite{goldman}), hence such a disk is often called a
Giraud disk (see~\cite{deraux4445}).

The existence of a third bisector containing $\B_1\cap\B_2$ may sound
mysterious at first, but it follows at once from the coequidistance
condition. Indeed, let $x_0$ be the intersecton point of the complex
spines $\Sigma_1$ and $\Sigma_2$, and let $x_j$, $j=1,2$ denote its
reflection across the real spine $\sigma_j$. Then $\B_j=\B(x_0,x_j)$,
and clearly $\B_1\cap \B_2$ is contained in $\B(x_1,x_2)$. The content
of Giraud's theorem is that these three bisectors are the only ones
containing $\B_1\cap \B_2$.

If the complex spines do not intersect, then they have a unique common
perpendicular complex line $\mathcal{T}$. This complex line is a slice
of $\B_1$ if and only if the real spine of $\Sigma_1$ goes through
$\Sigma_1\cap\mathcal{T}$ (and similarly for the real spine of
$\B_2$). This gives a simple criterion to check whether bisectors with
ultraparallel complex spines have a complex slice in common (this
happens if the extended real spines intersect). When this happens, the
bisectors are called {\bf cotranchal}. One should beware that when
this happens, the intersection can be strictly larger than the common
slice (but there can be at most one complex slice in common).

The slice parameters above allow an easy parametrization of the
intersection of the extors containing the bisectors, provided the
bisectors do not share a slice, which we now assume (this is enough
for the purposes of the present paper). In this case, the intersection
in projective space can be parametrized in a natural way by the
Clifford torus $S^1\times S^1\subset \C^2$. Specifically $(z_1,z_2)\in
S^1\times S^1$ parametrizes the vector orthogonal to
$\overline{z}_1\tilde{p}_0-\tilde{p}_1$ and
$\overline{z}_2\tilde{p}_2-\tilde{p}_3$. This vector can be written as
$$
(\overline{z}_1\tilde{p}_0-\tilde{p}_1)\boxtimes(\overline{z}_2\tilde{p}_2-\tilde{p}_3)
$$ 
in terms of the Hermitian box product, see~p. 43
of~\cite{goldman}. This can be rewritten in the form
\begin{equation}\label{eq:paramTorus}
V(\alpha,\beta)=c_{13} + z_1 c_{31} + z_2 c_{21} + z_1z_2 c_{02}
\end{equation}
where $c_{jk}$ denotes $p_j\boxtimes p_k$.

The intersection of the bisectors (rather than the extors) is given by
solving the inequality
$$
  \langle V(z_1,z_2),V(z_1,z_2)\rangle<0.
$$ 
The corresponding equation $\langle V(z_1,z_2),V(z_1,z_2)\rangle=0$ is
quadratic in each variable. It is known (see the analysis
in~\cite{goldman}) that the intersection has at most two connected
components. This becomes a bit simpler in the coequidistant case (then
one can take $p_0=p_2$, so that $c_{02}=0$), where the equation is
actually quadratic, rather than just quadratic in each variable.

Note that the intersection of three bisectors also has a simple
implicit parametrization, namely the intersection of $\B_1\cap\B_2$
with a third bisector $\B(q_1,q_2)$ has an equation
\begin{equation}\label{eq:raw}
|\langle V(z_1,z_2), \tilde{q}_1\rangle |^2 = 
   |\langle V(z_1,z_2), \tilde{q}_2\rangle |^2
\end{equation}
where $\tilde{q}_j$ are lifts of $q_j$ with the same square norm.

This implicit equation can be used to obtain piecewise
parametrizations for the corresponding curves, using either $z_1$ or
$z_2$ as a parameter. This is explained in detail in~\cite{deraux4445},
we briefly review some of this material.

Note that $\langle V(z_1,z_2), \tilde{q}_1\rangle$ is affine in each
variable (in the coequidistant case it is even affine in
$(z_1,z_2)$). This means that for a given $z_1$ with $|z_1|=1$,
finding the corresponding values of $z_2$ amounts to finding the
intersection of two Euclidean circles. Specifically, the equation has
the form
$$
|a_0(z_1)+a_1(z_1)z_2|^2=|b_0(z_1)+b_1(z_1)z_2|^2,
$$
which can be rewritten as
$$
2\Re( (\overline{a}_0a_1-\overline{b}_0b_1) z_2 ) = |a_0|^2+|a_1|^2-|b_0|^2-|b_1|^2,
$$
or simply in the form
\begin{equation}\label{eq:line-circle}
  \Re(\mu z_2)=\nu.
\end{equation}
Using the fact that $|z_1|=1$, we can write $\mu=\mu(z_1)$ and
$\nu=\nu(z_1)$ as affine functions in $(z_1,\overline{z}_1)$.

It follows from elementary Euclidean geometry (simply intersect the
circle of radius $|\mu|$ centered at the origin with the line
$\Re(z)=\nu$) that equation~\eqref{eq:line-circle} has a solution
$z_2$ with $|z_2|=1$ if and only if
\begin{equation}\label{eq:munu}
  |\mu|^2\geq \nu^2.
\end{equation}
If there is a $z_1$ such that $\mu=\nu=0$, then $z_2$ can of course be
chosen to be arbitrary (this happens when two of the three bisectors
share a slice). Otherwise, there is a single value of $z_2$
satisfying~\eqref{eq:line-circle} if and only if equality holds
in~\eqref{eq:munu}.

Of course the inequality $|\mu|^2\geq \nu^2$ can also be reinterpreted in
terms of the sign of the discriminant of a quadratic equation, since when $[z|=1$,
  $\mu z + \overline{\mu}\overline{z}=2\nu$ is equivalent to 
$$
  \mu z^2 - 2\nu z +\overline{\mu}=0.
$$

The determination of the projection of the curve~\eqref{eq:raw} onto
the $z_1$-axis of the Giraud torus amounts to the determination the
values of $z_1$, $|z_1|=1$ where there exists a $z_2$
satisfying~\eqref{eq:raw} and $|z_2|=1$. According to the previous
discussion, this amounts to finding where equality holds
in~\eqref{eq:munu}, which yields a polynomial equation in $z_1$. This
can be somewhat complicated, especially because polynomials can have
multiple roots.

On the intervals of the argument of $z_1$ corresponding to the
projection onto the $z_1$-axis of the curve defined by~\eqref{eq:raw}
(we remove the points where $\mu=\nu=0$ is arbitrary), we obtain a
nice piecewise parametrization for the curve, namely
\begin{equation}\label{eq:param}
  z_2 = \frac{\nu\pm i \sqrt{|\mu|^2-\nu^2}}{\mu}.
\end{equation}
This equation is problematic for numerical computations mainly when
$|\mu|$ is close to $\nu$. In that case, one can switch variables and
use $z_2$ rather than $z_1$ as the parameter.

All the above computations are fairly simple, but some care is needed
when performing them in floating point arithmetic. The main point that
allows us to perform somewhat sophisticated computations in our proofs
is the polynomial character of all equations, and the following.
\begin{prop}
  The group $\Gamma$ consists of matrices in $\textrm{GL}_3(K)$, where
  $K=Q(i\sqrt{7})$.
\end{prop}
Our fundamental domain is defined based on fixed points of certain
elliptic or parabolic elements in the group, whose coordinates can be
chosen to lie in $K$, so we will be able to choose the coefficients of
all the above polynomial parametrizations to lie in $K$. This allows
us to compute all relevant quantities to arbitrary precision; we will
treat some explicit sample computations in an appendix
(section~\ref{sec:appendix}).

Note that when the solution set of an equation of the
form~\eqref{eq:raw} is non empty, its dimension could in general be
$0,1$ or $2$.  Giraud's theorem (see Theorem~\ref{thm:giraud}) gives a
fairly general characterization of which bisectors can give a set of
dimension 2.

In the bisector intersections that appear in the present paper, we
will encounter situations where the solution set of~\eqref{eq:raw} is
a curve in the Clifford torus, but that intersects the closure in
$\chb 2$ of the Giraud disk only in a point at infinity. Among other
situations, this happens when the spinal spheres at infinity of
certain pairs of bisectors are tangent.

Clearly floating point arithmetic will give absolutely no insight
about such situations, so we will use geometric arguments instead. An
important geometric argument is the following result, proved by
Phillips in~\cite{phillips}:
\begin{prop}\label{prop:unipotent}
  Let $A$ be a unipotent isometry, and let $p_0\in \ch 2$. Then
  $\B(p_0,Ap_0)\cap \B(p_0,A^{-1}p_0)$ is empty.  The extension to
  $\partial_\infty\ch 2$ of these bisectors intersect precisely in the
  fixed point of $A$, in other words the spinal spheres for the above
  two bisectors are tangent at that fixed point.
\end{prop}
As we will see in the appendix (section~\ref{sec:appendix}), Phillips'
result allows to take care of most, but not all tangencies.

\section{Boundary unipotent representations} \label{sec:threereps}

We recall part of the results from~\cite{falbelFigure8}, using the
notation and terminology from section~\ref{sec:intro}, so that $M$
denotes the figure eight knot complement. We will interchangeably use
the following two presentations for $\pi_1(M)$:
\begin{equation}\label{eq:presfig8}
\langle\ g_1, g_2,
g_3\ |\ g_2=[g_3,g_1^{-1}],\ g_1g_2=g_2g_3\ \rangle
\end{equation}
and
$$
\langle\ a,b,t\ |\ tat^{-1}=aba,\ tbt^{-1}=ab\ \rangle. 
$$ 
The second presentation can be obtained from the first one by
setting $a=g_2$, $b=[g_2,g_3^{-1}]$ and $t=g_3$. Note that $a$ and $b$
generate a free group $F_2$, and the second presentation exhibits
$\pi_1(M)$ as the mapping torus of a pseudo-Anosov element of the
mapping class group of $F_2$; this comes from the fact that the figure
eight knot complement fibers over the circle, with once
punctured tori as fibers.

Representatives of the three conjugacy classes of representations of
$\pi_1(M)$ with unipotent boundary holonomy are the following
(see~\cite{falbelFigure8} pages 102-105).  We only give the image of
$g_1$ and $g_3$, since they clearly generate the group.
$$
 \rho_1(g_1)=\begin{pmatrix}
           1&  1 &     -\frac{1}{2}- \frac{\sqrt{3}}{2}i\\
 0 & 1 & -1  \\
0 & 0 & 1\\
           \end{pmatrix},
\  \
           \rho_1( g_3)=\begin{pmatrix}
           1&  0 &   0\\
 1 & 1 & 0  \\
  -\frac{1}{2}- \frac{\sqrt{3}}{2}i & -1 & 1\\
           \end{pmatrix}.
           $$
 $$
 \rho_2(g_1)=\begin{pmatrix}
           1&  1 &   -\frac{1}{2}- \frac{\sqrt{7}}{2}i\\
 0 & 1 & -1  \\
0 & 0 & 1\\
           \end{pmatrix},
 \ \
            \rho_2(g_3)=\begin{pmatrix}
           1&  0 &   0\\
 -1 & 1 & 0  \\
-\frac{1}{2}+ \frac{\sqrt{7}}{2}i& 1 & 1\\
           \end{pmatrix}.
           $$
$$
\rho_3( g_1)=\begin{pmatrix}
           1&  1 &   -1/2\\
 0 & 1 & -1  \\
0 & 0 & 1\\
           \end{pmatrix},
\ \             \rho_3(g_3)=\begin{pmatrix}
           1&  0 &   0\\
\frac{5}{4}- \frac{\sqrt{7}}{4}i  & 1 & 0  \\
-1& -\frac{5}{4}- \frac{\sqrt{7}}{4}i & 1\\
           \end{pmatrix}.
           $$
           
For completeness, we state the following result (the main part of
which was already proved in~\cite{falbelFigure8}).
\begin{prop}\label{prop:unique}
  For any irreducible representation $\rho:\pi_1(M)\rightarrow \pu21$
  with unipotent boundary holonomy, $\rho$ (or $\overline{\rho}$) is
    conjugate to $\rho_1$, $\rho_2$ or $\rho_3$.
\end{prop}
\pf \ We follow the beginning of section~5.4 in \cite{falbelFigure8}.
To prove this statement, we mainly need to complete the argument there
to exclude non generic cases.

Let $\rho$ be as in the statement of the proposition. In order to
avoid cumbersome notation, we will use the same notation as in the
introduction for the image of $g_1$, $g_2$ and $g_3$ under $\rho$, and
write $G_k=\rho(g_k)$.
 
We first observe that one of the boundary holonomy generators is given
by ${g_1}^{-1}{g_2}={g_1}^{-1}g_3{g_1}^{-1}{g_3}^{-1}g_1$.  This is
conjugate to ${g_1}^{-1}$ so $G_1=\rho(g_1)$ is unipotent by
assumption. Moreover, $g_1$ is conjugate to $g_3$, which implies that
$G_3=\rho(g_3)$ is unipotent as well.
  
Let $p_1$ and $p_2$ be the parabolic fixed points of $G_1 =\rho(g_1)$
and $G_3 =\rho(g_3)$, respectively.  We may assume that $p_1\neq p_2$
otherwise the representation would be elementary (hence not
irreducible).
  
Define $q_1 = {G_1}^{-1}(p_2)$ and $q_3=G_3(p_1)$.  By Lemma~5.3 in
\cite{falbelFigure8} (which uses only the presentation for $\pi_1(M)$,
see~\eqref{eq:presfig8}), 
 $$
      {G_3} {G_1}^{-1}(p_2 ) = {G_1}^{-1}G_3 (p_1).
 $$
We define $q_2$ as the point on both sides of the above equality.
  
If $p_1, p_2, q_1, q_2$ and $p_1, p_2, q_2, q_3$ are in general
position (that is, no three points belong to the same complex line)
these quadruples are indeed parametrized by the coordinates
from~\cite{falbelFigure8}, and these coordinates must be solutions of
the compatibility equations, so $\rho$ must be conjugate to some
$\rho_j$ (or its complex conjugate).
  
If the points are not in general position we analyze the
representation case by case.
    
The first case is when $q_1={G_1}^{-1}(p_2)$ belongs to the boundary
of the complex line through $p_1$ and $p_2$.  Without loss of
generality, we may assume $p_1=\infty$ and $p_2=(0,0)$ in Heisenberg
coordinates.  As $G_1$ preserves the complex line between $p_1$ and
$p_2$ it has the following form:
    $$
    G_1=\begin{pmatrix}
           1&  0 &   \frac{it}{2}\\
           0 & 1 & 0  \\
           0& 0 & 1\\
           \end{pmatrix}.
           $$
    We then write 
   $$
    G_3=\begin{pmatrix}
           1&  0 &  0\\
           z & 1 & 0 \\
           - \frac{|z|^2}{2}+\frac{is}{2}& -\bar z & 1\\
          \end{pmatrix}.
  $$ 
with $z\neq 0$ (otherwise the representation would be reducible).
    Now, the equation $${G_3}^{-1} {G_1}^{-1}(p_2 ) = {G_1}^{-1}G_3
    (p_1 )$$ gives
   $$
   \begin{pmatrix}
         -\frac{it}{2} \\
         -\frac{izt}{2}\\
         \frac{it|z|^2}{4}+\frac{ts}{4}+1\\
           \end{pmatrix}= 
           \lambda 
            \begin{pmatrix}
       \frac{it|z|^2}{4}+\frac{ts}{4}+1 \\
 z\\
 -\frac{|z|^2}{2}+\frac{is}{4}\\
           \end{pmatrix} .
  $$ 
One easily checks that this equation has no solutions with $z\neq 0$.
Therefore $q_1$ is not in the complex line defined by $p_1$ and $p_2$.

Analogously, $q_3=G_3(p_1)$ cannot be in that complex line
either. Now, from the gluing pattern in Figure~\ref{fig:tetrahedra},
we obtain that $p_1,q_1, q_2$ and $p_2,q_2,q_3$ are in general
position.  It remains to verify that $p_2,q_1,q_2$ are in general
position.  We write
  $$
  (p_2,q_1,q_2)=(p_2,{G_1}^{-1}(p_2 ),{G_1}^{-1}G_3 (p_1 ) )={G_1}^{-1}G_3({G_3}^{-1}G_1(p_2),p_2,p_1)
  $$
 But if $
 ({G_3}^{-1}G_1(p_2),p_2,p_1)
 $
are on the same complex line then, again, we obtain equations which
 force $p_1,p_2, q_1$ to be in the same line.  \cqfd

In fact it is not hard to show that there are no reducible
representations apart from elementary ones (still assuming the boundary
holonomy to be unipotent). The relator relation then implies that these
elementary representations must satisfy $\rho(g_1)=\rho(g_3)$, hence
the image of the representation is in fact a cyclic group.

\section{A Dirichlet domain for $\Gamma$} \label{sec:statement}

From this point on, we mainly focus on the representation $\rho_2$
(see the discussion in the introduction, and section~\ref{sec:mcg}).
We write $\Gamma=\Gamma_2$ and
$$
G_1=\rho_2(g_1),\quad G_2=\rho_2(g_2),\quad G_3=\rho_2(g_3).
$$

The combinatorics of Dirichlet domains depend significantly on their
center $p_0$, and there is of course no canonical way to choose this
center. We will choose a center that produces a Dirichlet domain with
very few faces, and that has a lot of symmetry (see
section~\ref{sec:symmetry}), namely the fixed point of $G_2$.

Recall that $G_2=[G_3,G_1^{-1}]$, and this can easily be computed to
be
$$
G_2=\left(
  \begin{array}{ccc}
    2 & \frac{3}{2}-i\frac{\sqrt{7}}{2} & -1 \\
    -\frac{3}{2}-i\frac{\sqrt{7}}{2} & -1 & 0 \\
    -1 & 0 & 0 \\
  \end{array}
\right)
$$ 
It is easy to check that $G_2$ is a regular elliptic element of
order $4$, whose isolated fixed point is given in homogeneous
coordinates by
$$
\tilde{p}_0=(1,-(3+i\sqrt{7})/4,-1).
$$ 
Note that no nontrivial power of $G_2$ fixes any point in
$\partial_\infty \ch 2$ ($G_2$ and $G_2^{-1}$ are regular elliptic,
and $G_2^2$ is a complex reflection in a point).

Recall from section~\ref{sec:dirichlet} that, for any subset
$S\subset\Gamma$, $E_S$ denotes the Dirichlet domain centered at
$p_0$; the faces of $E_S$ are given by intersections of the form
$$
E_S\cap \B(p_0,\gamma p_0)
$$ 
that have non empty interior in $\B(p_0,\gamma p_0)$ (we refer to
such a face as being associated to the element $\gamma$).

As a special case, $E_\Gamma$ denotes the Dirichlet domain for
$\Gamma$ centered at $p_0$, and $E_S$ denotes an a priori larger
domain taking into account only the faces coming from $S$ rather than
all of $\Gamma$.

From this point on, we will always fix the set $S$ to be the following
set of eight group elements:
\begin{equation}\label{eq:S}
S = \{G_2^kG_1G_2^{-k}, G_2^kG_3^{-1}G_2^{-k}, k=0,1,2,3\}.
\end{equation}
Since for the remainder of the paper we will always use the same set
$S$, we simply write
$$
E = E_S.
$$ 
Note also that it follows from simple relations in the group that $S$
is a \emph{symmetric} generating set (in the sense that it is closed
under the operation of taking inverses in the group), even though this
may not be obvious from the above description. For now we simply refer
to the second column of Table~\ref{tab:faces}, where the relevant
relations in the group are listed.

With this notation, what we intend to prove is the following (which
will be key to the proof of Theorem~\ref{thm:main}).
\begin{thm}\label{thm:submain}
  The Dirichlet domain $E_\Gamma$ centered at $p_0$ is equal to
  $E$. In particular, $E_\Gamma$ has precisely eight faces, namely
  the faces of $E_\Gamma$ associated to the elements of
  $S$, which are listed in~\eqref{eq:S}.
\end{thm}
As outlined in section~\ref{sec:dirichlet}, in order to prove that
$E_\Gamma=E$, we will start by determining the precise combinatorics
of $E$, then apply the Poincar\'e polyhedron theorem in order to prove
that $E$ is a fundamental domain for $\Gamma$ modulo the action of the
finite group $H$.

Note that $E$ is indeed not a fundamental domain for $\Gamma$, since
by construction it has a nontrivial stabilizer (powers of $G_2$ fix
the center of $E$, hence they must preserve $E$). It is a fundamental
domain for the coset decomposition of $\Gamma$ into left cosets of the
group $H$ of order $4$ generated by $G_2$ (see
section~\ref{sec:dirichlet}), and this suffices to produce a
presentation for $\Gamma$, see section~\ref{sec:presentation}.  One
can deduce from $E$ a fundamental domain for $\Gamma$, by taking
$E\cap F$ where $F$ is any fundamental domain for $H$. We omit the
details of that construction, since they will not be needed in what
follows.

\begin{dfn}\label{def:boundingBisectors}
We write $\B_1$,\dots,$\B_8$ for the bisectors bounding $E$, numbered
as in Table~\ref{tab:faces}. For each $k$, we denote by $\OB_k$ the
closure of $\B_k$ in $\chb 2=\ch 2\cup\partial_\infty \ch2$. We write
$b_k$ for the intersection $\B_k\cap E$, and $\overline{b}_k$ for the
closure of that face in $\chb 2$.
\end{dfn}
We will sometimes refer to the bisectors $\B_j$ as the {\bf bounding
  bisectors}.

\begin{table}[htbp]
\centering
  \begin{tabular}{c|c|c|c}
  Element of $S$    & Bisector              &  Face   & Vertices\\
  \hline
   $G_1$                                & $\B_1 $              &  $b_1$  &  $p_1,p_2,q_3,q_4$\\
   $G_3^{-1}$                           & $\B_2$               &  $b_2$ & $p_2,q_4,q_1,p_1$\\
   $G_2G_1G_2^{-1}$                     & $G_2\B_1=\B_3$       &  $b_3$ & $p_4,p_1,q_4,q_1$\\ 
   $G_2G_3^{-1}G_2^{-1}=G_1^{-1}$       & $G_2\B_2=\B_4 $      &  $b_4$ & $p_1,q_1,q_2,p_4$\\
   $G_2^{2}G_1G_2^{2}=G_2^{-1}G_3G_2$   & $G_2^2\B_1=\B_5   $  &  $b_5$ & $p_3,p_4,q_1,q_2$  \\
   $G_2^2G_3^{-1}G_2^{2}=G_2G_1^{-1}G_2^{-1}$ & $G_2^2\B_2=\B_6 $  &  $b_6$ & $p_4,q_2,q_3,p_3$  \\
   $G_2^{-1}G_1G_2=G_3$                 & $G_2^{-1}\B_1=\B_7   $  &  $b_7$  & $p_2,p_3,q_2,q_3$  \\
   $G_2^{-1}G_3^{-1}G_2$                & $G_2^{-1}\B_2=\B_8 $ &  $b_8$ & $p_3,q_3,q_4,p_2$
  \end{tabular}
  \caption{Notation for the eight faces of the Dirichlet domain; the
  face associated to an element $\gamma\in S$ is contained in
  $\B(p_0,\gamma p_0)$, see section~\ref{sec:dirichlet}.
  The equalities in the first column follow from the relation
  $G_1G_2=G_2G_3$. 
  The notation for vertices will be explained in
  section~\ref{sec:vertices}.} \label{tab:faces}
\end{table}

\subsection{Symmetry} \label{sec:symmetry}

Note that $S$ is by construction invariant under conjugation by $G_2$,
which fixes $p_0$, so $E$ is of course $G_2$-invariant. In
particular, it has at most $2$ isometry types of faces; in fact all
its faces are isometric, as can be seen using the involution
$$
I=\left(\begin{matrix}
  0 & 0 & 1\\
  0 & -1 & 0\\
  1 & 0 & 0
\end{matrix}\right).
$$ 
This is not an element of $\Gamma$, but it can easily be
checked that it normalizes $\Gamma$ by using the conjugacy
information given in Proposition~\ref{prop:conjugacy}.
\begin{prop}\label{prop:conjugacy}
  $$IG_1I=G_3^{-1}$$
  $$IG_2I=G_2^{-1}$$
\end{prop}
This proposition shows that the group generated by $I$ and $G_2$ has
order $8$, and this group of order $8$ stabilizes $E$ (the formula
given above for $p_0$ makes it clear that it is fixed by
$I$). Finally, note that Proposition~\ref{prop:conjugacy} makes it
clear that $I$ exchanges the faces $b_1$ and $b_2$ (see
Table~\ref{tab:faces} for notation).

\subsection{Vertices of $E$} \label{sec:vertices}

In this section we describe certain fixed points of unipotent elements
in the group, which will turn out to give the list of all vertices of
$E$ (this claim will be justified in the end of
section~\ref{sec:combinatorics}, see
Proposition~\ref{prop:vertices}). We use the numbering of faces (as
well as bisectors that contain these faces) given in
Table~\ref{tab:faces}. We start mentioning that $G_1$ clearly maps
$\B_4=\B(p_0,G_1^{-1}p_0)$ to $\B_1=\B(p_0,G_1p_0)$. Since $G_1$ is
unipotent, Proposition~\ref{prop:unipotent} shows that the
corresponding bisectors have empty intersection, and their spinal
spheres are tangent at the fixed point of $G_1$.

The latter is clearly given by $$p_1=(1,0,0),$$ and it is easy to
check that this point is on the closure of precisely four of the
bisectors that bound the Dirichlet domain, namely $\OB_1$, $\OB_2$,
$\OB_3$ and $\OB_4$. The fact that it is in $\OB_1$ and $\OB_4$ is
obvious, the other ones can be checked by explicit
computation. Indeed, we have
$$
G_3\ p_1 = (1,-1,\frac{-1+i\sqrt{7}}{2}),
$$
$$
G_1^{-1}G_2^{-1}\ p_1= (\frac{1-i\sqrt{7}}{2},-1,-1)
$$
hence
$$
|\langle p_1,G_3^{-1}p_0\rangle| = |\langle G_3\ p_1,p_0\rangle| = 1 = |\langle p_1,p_0\rangle|
$$
$$
|\langle p_1,G_2G_1G_2^{-1}p_0\rangle| = |\langle G_1^{-1}G_2^{-1}p_1,p_0\rangle| = 1 = |\langle p_1,p_0\rangle|.
$$

Similarly, the bisectors $\B_2$ and $\B_5$ have tangent spinal
spheres, and this comes from the fact that $G_2^{-1}G_3$ is unipotent
(which can be checked by direct calculation). Indeed, this isometry
sends $\B_2=\B(p_0,G_3^{-1}p_0)$ to
$\B(G_2^{-1}G_3p_0,G_2^{-1}p_0)=\B(G_2^{-1}G_3G_2p_0,p_0)=\B_5$.

We call $q_1$ the fixed point of $G_2^{-1}G_3$, which can easily be
computed to be given by
$$
q_1=(\frac{-1+i\sqrt{7}}{2},1,1).
$$ 
One verifies directly that this point is on the closure of
precisely four bounding bisectors, namely $\OB_2$, $\OB_3$,
$\OB_4$ and $\OB_5$.

Now applying $G_2$ to $p_1$ and $q_1$, we get eight specific fixed
points of unipotent elements in the group which are all tangency
points of certain spinal spheres. We define points $p_k$, $q_k$ for
$k=1,\dots,4$ by
$$
p_{k}=G_2 p_{k+1}; \quad q_{k+1}=G_2 q_k.
$$ 
Beware that $G_2$ raises the indices of $q$-vertices, whereas it
lowers the indices of the $p$-vertices; this somewhat strange
convention is used for coherence with the notation
in~\cite{falbelwang}.

Perhaps surprisingly, the eight tangency points will turn out to give
all the vertices of the Dirichlet domain. We summarize the results in
the following.
\begin{prop}\label{prop:tangency}
  There are precisely eight pairs of tangent spinal spheres among the
  boundary at infinity of the bisectors bounding the Dirichlet
  domain. The list of points of tangency is given in
  Table~\ref{tab:vertexLabels}.
\begin{table}[htbp]
\centering
  \begin{tabular}{c|c|c|c}
    Vertex &  Fixed by & tangent spinal spheres & Other faces\\\hline 
    $p_1$ & $G_1$ & $\B_1$,\ $\B_4$ & $\B_2$,\ $\B_3$\\
    $p_2$ & $G_3$ & $\B_7$,\ $\B_2$ & $\B_8$,\ $\B_1$\\
    $p_3$ & $G_2^{-1}G_3G_2$ & $\B_5$,\ $\B_8$ & $\B_6$,\ $\B_7$\\
    $p_4$ & $G_2G_1G_2^{-1}$ & $\B_3$,\ $\B_6$ & $\B_4$,\ $\B_5$\\
    $q_1$ & $G_3^{-1}G_2$ & $\B_2$,\ $\B_5$ & $\B_3$,\ $\B_4$\\
    $q_2$ & $G_1^{-1}G_2$ & $\B_4$,\ $\B_7$ & $\B_5$,\ $\B_6$\\
    $q_3$ & $G_2G_1^{-1}$ & $\B_6$,\ $\B_1$ & $\B_7$,\ $\B_8$\\
    $q_4$ & $G_3G_1$ & $\B_8$,\ $\B_3$ & $\B_1$,\ $\B_2$
  \end{tabular}
\caption{The vertices of $E$ at infinity, given by a unipotent element
  that fixes them. See also the list of vertices that lie on each face
  given in Table~\ref{tab:faces}.}\label{tab:vertexLabels}
\end{table}
\end{prop}

\begin{pf}
The claim about tangency has already been proved, we only justify the
fact that the points in the $G_2$-orbit of $p_1$ and $q_1$ are indeed
stabilized by the unipotent element given in
Table~\ref{tab:vertexLabels}. This amounts to checking that the
unipotent elements claimed to fix the points $p_j$ (resp. those
claimed to fix the points $q_j$) are indeed conjugates of each other
under powers of $G_2$.

This can easily be seen from the presentation of the group (in fact
the relations $G_1G_2=G_2G_3$ and $(G_1G_2)^3=1$ suffice to check
this). For instance, $G_2(p_4)=p_3$ because, using standard word
notation in the generators where $1=G_1$,
$\overline{1}=G_1^{-1}$, we have
$$
2\cdot 2\bar1\,\bar2 \cdot \bar2 
= \bar 2 \bar 2 \cdot \bar 1 \cdot 2 2
= \bar2 \cdot \bar2 \bar 1 \cdot 2 2 
= \bar2 \cdot  \bar3 \bar 2 \cdot 2 2 
= \bar2\bar3 2.
$$

Similarly, $G_2(q_3)=q_4$ because
$$
2\cdot 2\bar1\cdot\bar2 = 2 2 \cdot \bar 1\bar 2 = 2 2 \cdot 2 1 2 1 = \bar2 \cdot 12 \cdot 1 = \bar2 \cdot 2 3 \cdot 1 = 3 1.
$$
The other conjugacy relations are handled in a similar fashion.
\end{pf}

\subsection{Combinatorics of $E$} \label{sec:combinatorics}

We now go into the detailed study of the combinatorics of $E$. 

\begin{figure}[htbp]
  \centering
  \includegraphics[width=0.8\textwidth]{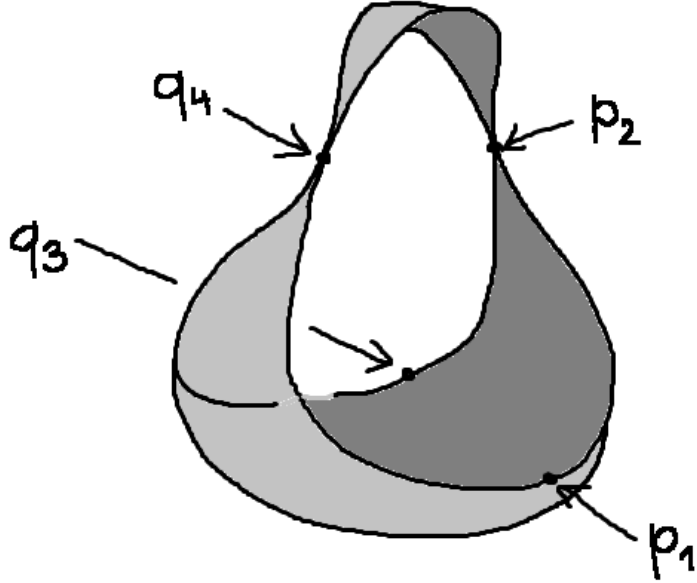}
  \caption{A schematic picture of $\partial_\infty b_1$. The face
    $\overline{b}_1$ also has two finite 2-faces, given by the Giraud
    disks $\OB_8\cap\OB_1$ and $\OB_1\cap \OB_2$ (only their boundary
    circle is draw in the picture).
    The face $\overline{b}_1$ has precisely four vertices, all
    in the ideal boundary (they are the fixed points of $G_1$, $G_3$,
    $G_3G_1$, and $G_1G_2^{-1}$, see Tables~\ref{tab:faces}
    and~\ref{tab:vertexLabels}).} 
  \label{fig:b}
\end{figure}

The results of section~\ref{sec:symmetry} show that it is enough to
determine the combinatorics of a single face of $E$, say
$b_1=E\cap\B_1$, and its incidence relation to all other faces.
\begin{prop}\label{prop:comb-b1}
  The closure $\overline{b}_1$ of $b_1$ in $\chb 2$ has precisely
  three 2-faces, two finite ones and one on the spinal sphere
  $\partial_\infty \B_1$.
  \begin{enumerate}
    \item The finite 2-faces are the given by the (closure of the)
      Giraud disks $\OB_1\cap \OB_2$, $\OB_1\cap \OB_8$;
    \item The 2-face on the spinal sphere $\partial_\infty \B_1$ is an
      annulus, pinched at two pairs of points on its boundary. The
      pinch points correspond to the fixed points of $G_3$ and
      $G_3G_1$.
  \end{enumerate}
  In particular, $\overline{b}_1$ intersects all faces $\OB_k$, $k\neq
  2,7$ in lower-dimensional faces.
\end{prop}
A schematic picture of the combinatorics of $\overline{b}_1$ is given
in Figure~\ref{fig:b}, where the shaded region corresponds to the
$2$-face of $\overline{b}_1$ at infinity (part~(2) of the
Proposition).
The Giraud disks mentioned in part~(1) of the Proposition intersect
only in two points in $\partial_\infty\ch 2$, not inside $\ch 2$ (see
Proposition~\ref{prop:vertices}).
The intersection pattern of the boundary at infinity of the eight
faces $\B_1,\dots,\B_8$ is somewhat intricate. Eight isometric copies
of the shaded region in Figure~\ref{fig:b} are glued according to the
pattern illustrated in Figure~\ref{fig:torusFig8}
(section~\ref{sec:combinatoricsAtInfinity}).

The general remark is that the claims in
Proposition~\ref{prop:comb-b1} can be proved using the techniques of
section~\ref{sec:bisectors}. In this section, we break up the proof of
Proposition~\ref{prop:comb-b1} into several lemmas
(Lemma~\ref{lem:bisectorIntersections},~\ref{lem:faceIntersections}),
and make these lemmas plausible by drawing pictures that can easily be
reproduced using the computer (and the parametrizations explained in
section~\ref{sec:bisectors}). The detailed proof will be postponed
until the appendix (section~\ref{sec:appendix}), since it relies on
somewhat delicate computations.

\begin{figure}[htbp]
  \centering
  \hfill 
  \subfigure[$\B_1\cap \B_2$]{\includegraphics[width=0.5\textwidth]{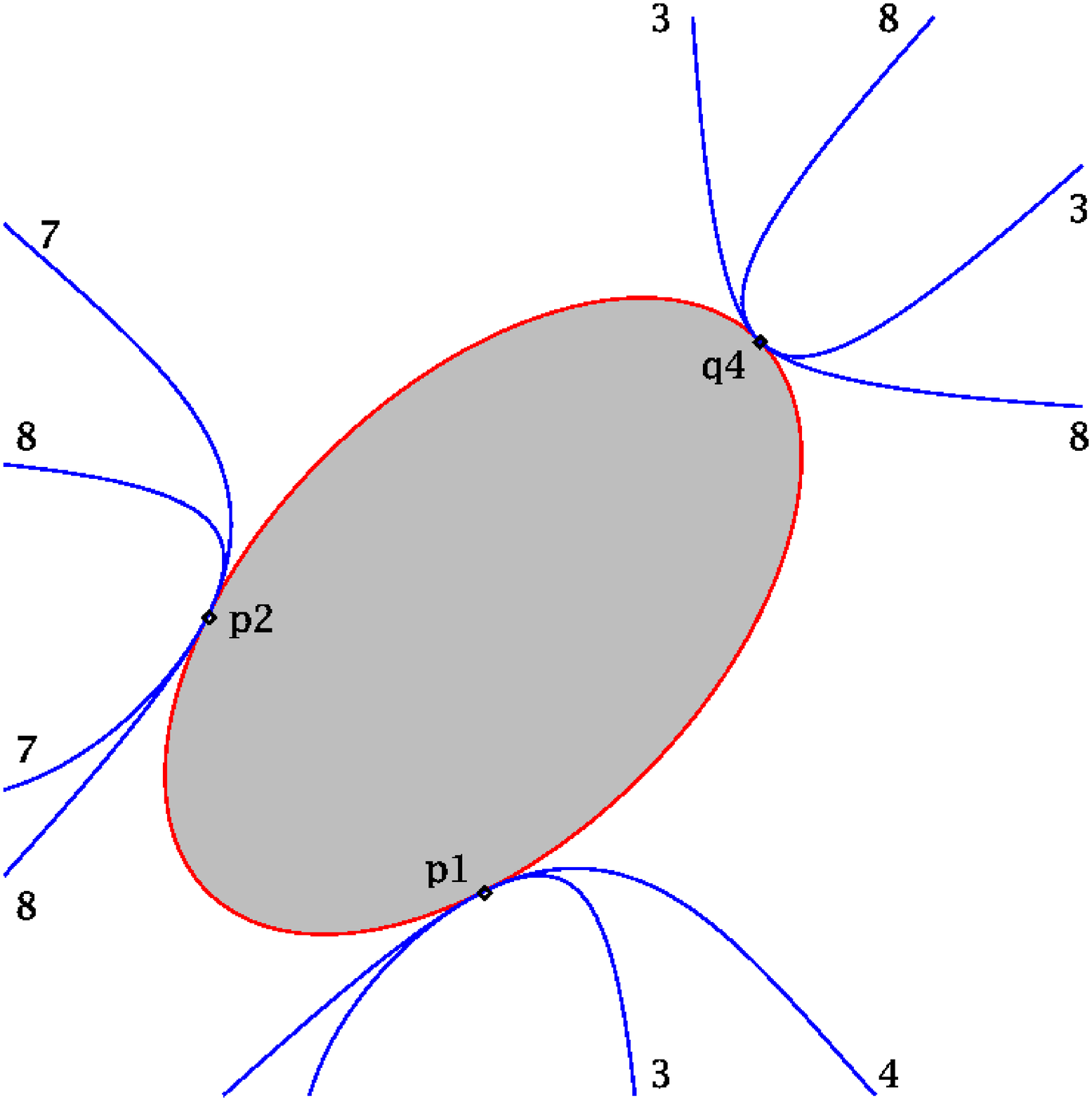}} \hfill
  \subfigure[$\B_1\cap \B_7$]{\includegraphics[width=0.45\textwidth]{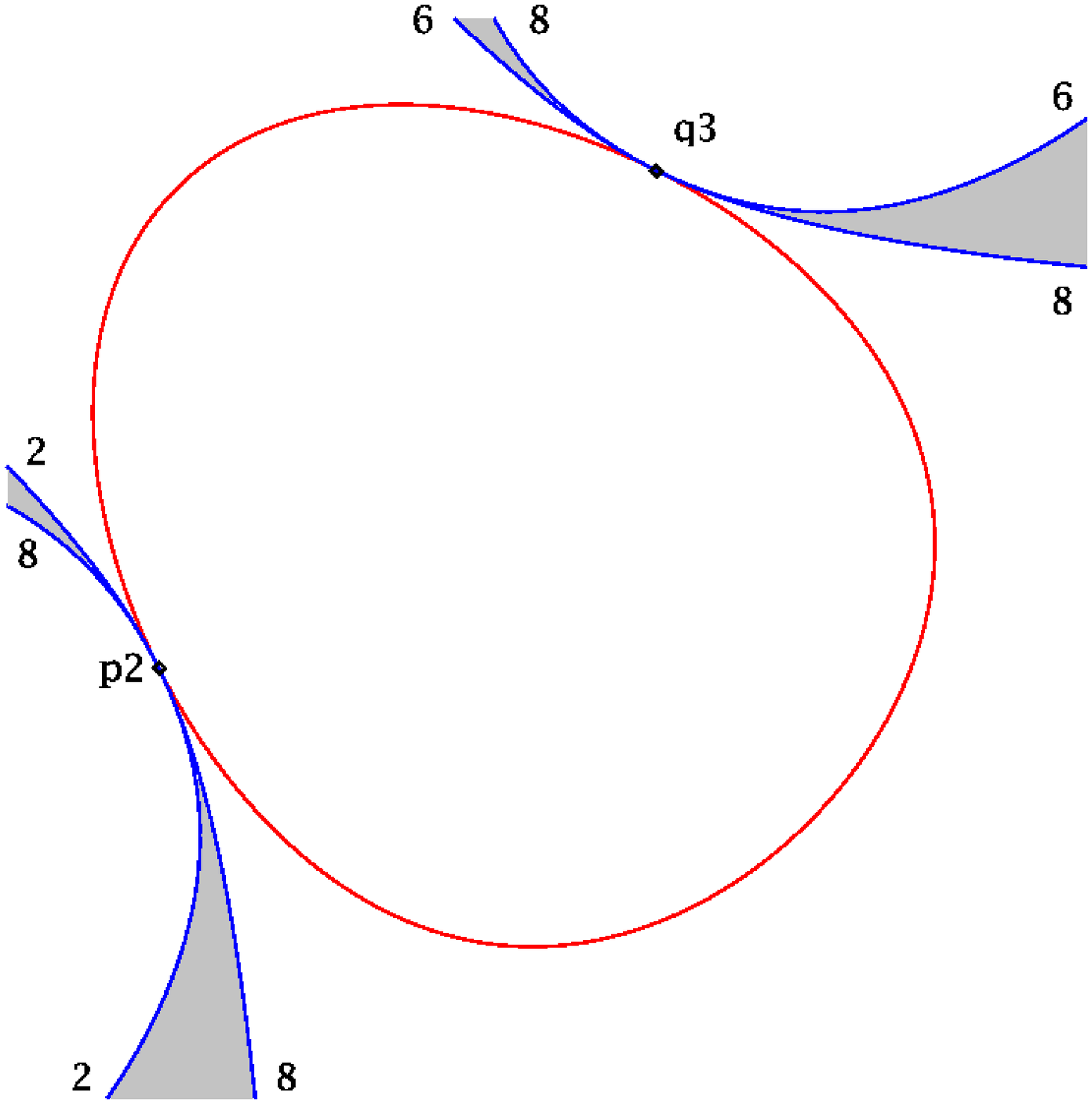}}\hfill
  \caption{Typical Giraud disk corresponding to the intersection of
    two bounding bisectors; the other curves are traces of the other
    $6$ bisectors.}\label{fig:disks}
\end{figure}

Since any two of the eight bisectors bounding $E$ are coequidistant,
their pairwise intersections are either empty, or diffeomorphic to a
disk (see section~\ref{sec:bisectors}). Recall that such disks are
either complex lines or Giraud
disks. Lemma~\ref{lem:bisectorIntersections} details the intersections
of $\B_1$ with the seven bisectors $\B_k$, $k\neq 1$. It can easily be
translated into a statement about $\B_2$ by using the involution $I$
(see section~\ref{sec:symmetry}), hence also about any $\B_j$ by using
powers of $G_2$.
\begin{lem}\label{lem:bisectorIntersections}
  $\B_1$ intersects exactly four of the seven other bisectors bounding
  $E$, namely $\B_7$, $\B_8$, $\B_2$ and $\B_3$. The corresponding
  intersections are Giraud disks.
\end{lem}
\begin{pf}
  The fact that $\B_1\cap\B_4$ and $\B_1\cap\B_6$ are empty follows
  from Proposition~\ref{prop:tangency}.  The fact that $\B_1\cap
  \B_5=\emptyset$ can be shown with direct computation, using the
  parametrization of the corresponding Giraud torus explained in
  section~\ref{sec:bisectors}. 
  The fact that the intersection of $\B_1$ with the four bisectors in
  the statement is indeed a Giraud disk can be done simply by
  exhibiting a point in that Giraud disk. Details will be given in
  the appendix (section~\ref{sec:boundingpairs}).
\end{pf}
The following statement is the analogue of
Proposition~\ref{lem:bisectorIntersections}, pertaining to face
(rather than bisector) intersections.
\begin{lem}\label{lem:faceIntersections}
  \begin{enumerate}
  \item $\B_1\cap \B_7\cap E$ and $\B_1\cap \B_3\cap E$ are empty.

  \item $\B_1\cap\B_2\cap E=\B_1\cap\B_2$ and $\B_1\cap \B_8\cap
  E=\B_1\cap \B_8$, and these are both Giraud disks.
  \end{enumerate}
\end{lem}
The proof of this statement will be given in the appendix
(section~\ref{sec:proof-faceint}). For now, we only show some pictures drawn
in spinal coordinates on the relevant Giraud disks, see
Figure~\ref{fig:disks}. 
For each of them we plot the trace on that Giraud disk of the other
six bisectors (see section~\ref{sec:bisectors} for a description of
how this can be done).  In the picture, we label each arc with the
index of the corresponding bisector (see the numbering in
Table~\ref{tab:faces}).

The fact that these pictures can indeed be trusted depends on the fact
that the curves have polynomial equations with entries in an explicit
number field $\Q(i\sqrt{7})$, as will be explained in detail in
section~\ref{sec:proof-faceint} of the appendix.

It follows from the previous analysis that the face $b_1$ has no
vertex in $\ch 2$, and that it has exactly four ideal vertices, or in
other words the closure $\overline{b}_1$ has four vertices. We
summarize this in the following proposition:
\begin{prop}\label{prop:vertices}
  $\overline{b}_1$ has precisely four vertices, all at infinity. They
  are given by $p_1$, $p_2$, $q_3$, $q_4$.
\end{prop}
\begin{pf}
  $p_2$ and $q_3$ are obtained as the only two points in the
  intersection $\OB_1\cap \OB_7$ (as before, bars denote
  the closures in $\ch 2\cup\partial_\infty \ch2$), see
  Figure~\ref{fig:disks}. Similarly, $p_1$ and $q_4$ are the two
  points in $\OB_1\cap \OB_3$.
\end{pf}
One can easily use symmetry to give the list of vertices of every
face. Each face has precisely four (ideal) vertices, see the last
column of Table~\ref{tab:faces}.

\section{The Poincar\'e polyhedron theorem for $E$} \label{sec:poincare}

This section is devoted to proving the hypotheses of the Poincar\'e
polyhedron theorem for the Dirichlet polyhedron $E$
(sections~\ref{sec:sidePairings},~\ref{sec:ridgeCycles}
and~\ref{sec:vertexCycles}), and to state some straightforward
applications (section~\ref{sec:presentation}).

\subsection{Side pairings}\label{sec:sidePairings}

We now check that opposite faces of $E$ (i.e. faces that correspond to
$\gamma$ and $\gamma^{-1}$, for $\gamma\in S$) are paired by the
isometry $\gamma$. It is enough to check this for $\gamma=G_1$, since
all others are obtained from this one by symmetry. More concretely, we
will check that $G_1$ maps $b_4$ to $b_1$, see Table~\ref{tab:faces}
for notation.

Recall that $\overline{b}_4$ has three facets, one on the ideal boundary
$\partial_\infty \ch 2$ and two given by the Giraud disks
$\B_4\cap\B_3$ and $\B_4\cap \B_5$.
\begin{prop}\label{prop:pairing}
  The isometry $G_1$ maps $\B_4\cap\B_3$ to $\B_1\cap\B_2$, and
  $\B_4\cap\B_5$ to $\B_1\cap\B_8$.
\end{prop}
\begin{pf}
  The Giraud disk $\B_4\cap\B_3$ is equidistant from $p_0$,
  $G_1^{-1}p_0$ and $G_2G_1G_2^{-1}p_0=G_2G_1p_0$, whereas
  $\B_1\cap\B_2$ is equidistant from $p_0$, $G_1p_0$ and $G_3^{-1}
  p_0$.

  Now $G_1(\B_4\cap\B_3)=\B_1\cap\B_2$ is equivalent to
  $$
  G_1G_2G_1 p_0 = G_3^{-1} p_0,
  $$ 
  which can easily be checked by direct computation. Equivalently,
  one may check that $G_3G_1G_2G_1 = G_2^2$.

  The fact that $G_1(\B_4\cap\B_5)=\B_1\cap\B_8$ follows similarly from 
  $$
  G_1G_2^{-1}G_3p_0 = G_2^{-1}G_3^{-1}p_0,
  $$
  or equivalently $G_3G_2G_1G_2^{-1}G_3=G_2^2$. 

  These relations in the group are of course easily obtained from the
  group presentation, but they can also be checked directly from the
  explicit matrices that appear in section~\ref{sec:threereps}.
\end{pf}

Proposition~\ref{prop:pairing} implies that $G_1$ maps $b_4$
isometrically to $b_1$. We will need more specific information about
the image of vertices under the side pairings (see the last column of
Table~\ref{tab:faces} for the list of vertices on each face, where the
quadruples of vertices are ordered in a consistent manner, i.e. the
side pairing maps the $j$-th vertex to the $j$-th vertex).
\begin{prop}
  The isometry $G_1$ maps the vertices of $b_4$ to vertices of face
  $b_1$. More specifically, $G_1(p_1)=p_1$, $G_1(p_4)=q_4$,
  $G_1(q_1)=p_2$ and $G_1(q_2)=q_3$.
\end{prop}
\begin{pf}
  The fact that $G_1(p_1)=p_1$ is obvious. The point $q_2=G_2(q_1)$ is
  the fixed point of $G_1^{-1}G_2$, so $G_1(q_2)$ is fixed by
  $G_2G_1^{-1}$, hence the latter point must be $q_3$ (see the second
  column in Table~\ref{tab:vertexLabels}).

  The fact that
  $G_1(q_1)=p_2$ follows from the fact that $g_2=[g_3,g_1^{-1}]$, since
  $$ 
  1 \cdot \bar{3} 2 \cdot \bar{1} = 1\bar3 \cdot [3,\bar1] \cdot \bar 1 = \bar{3}.
  $$

  Finally, the fact that $G_1(p_4)=q_4$ is equivalent to showing that
  $121\bar 2\bar 1$ and $31$ have the same fixed point. This follows
  from $(12)^2=(121)^3=id$ and $12=23$, since
  $$
  121\bar2\bar 1 = 121\cdot 1212 = (\bar1\bar2\bar1)^2 \cdot 1212 =\bar 1\bar 2\bar 1 2 = \bar 1\bar 3\bar 2 2 = \bar 1\bar 3.
  $$
\end{pf}

\subsection{Cycles of ridges}\label{sec:ridgeCycles}

It follows from Giraud's theorem (Theorem~\ref{thm:giraud}) that
the ridges of $E$ are on precisely three bisectors, hence there are
three copies of $E$ tiling its neighborhood. We only need to consider
the ridges $b_1\cap b_2$ and $b_1\cap b_8$, since the other ones are
all images of these two under the appropriate power of $G_2$.

The Giraud disk $\B_1\cap \B_2$ is equidistant from $p_0$, $G_1(p_0)$ and
$G_3^{-1}(p_0)$, and we apply $G_1^{-1}$ to this triple of points,
getting $G_1^{-1}p_0,p_0,G_1^{-1}G_3^{-1}p_0=G_2G_1p_0$, and bring it
back to $\B_1\cap \B_2$ by applying $G_2^{-1}$. This does not yield the
identity, but effects a cyclic permutation of the above three points:
$$
\begin{array}{c}
  p_0, G_1 p_0, G_3^{-1}p_0\\
  \downarrow G_2\\
  G^{-1} p_0, p_0, G_1^{-1}G_3^{-1} p_0\\
  \downarrow G_1\\
  G_3^{-1} p_0, p_0, G_1 p_0\\
\end{array}
$$ 
In other words, the corresponding cycle transformation is $G_1G_2$,
and the corresponding relation is
$$
(G_1G_2)^3=Id.
$$
Another geometric interpretation of this the following:
\begin{prop}\label{prop:tiling_b1_b2}
  A neighborhood of a generic point of $b_1\cap b_2$ is tiled by $E$,
  $G_1G_2(E)=G_1(E)$, $(G_1G_2)^{-1}(E)=G_3^{-1}(E)$.
\end{prop}

The Giraud disk $\B_1\cap\B_8$ is equidistant from $p_0$, $G_1(p_0)$
and $G_3 G_1(p_0)$. Again, we get an isometry in the group that
permutes these points cyclically:
$$
\begin{array}{c}
  p_0, G_1 p_0, G_3 G_1p_0\\
  \downarrow G_2^{2}\\
  G_1^{-1} p_0, p_0, G_1^{-1}G_3G_1 p_0\\ 
  \downarrow G_1\\
  G_3G_1 p_0, p_0, G_1 p_0
\end{array}
$$ 
which gives the relation
$$
(G_1G_2^2)^3=Id.
$$

The statement analogous to Proposition~\ref{prop:tiling_b1_b2} is the
following:
\begin{prop}\label{prop:tiling_b1_b8}
  A neighborhood of a generic point of $b_1\cap b_8$ is tiled by $E$,
  $G_2^2G_1(E)$ and $(G_2^2G_1)^{-1}(E)=G_1^{-1}(E)$.
\end{prop}

\subsection{Cycles of boundary vertices} \label{sec:vertexCycles}

As explained in section~\ref{sec:dirichlet}, we need to check that the
cycle transformations for all boundary vertices are parabolic. There
is only one cycle of vertices, since $G_2(q_k)=q_{k+1}$,
$G_2(p_k)=p_{k-1}$ (indices mod $4$), and we have
$$
G_3(p_1)=q_3.
$$ 

We check the geometry of the tiling of $\ch 2$ near $p_1$, which can
be deduced from the structure of ridges through that point (see
section~\ref{sec:ridgeCycles}).  Recall that $p_1$ is on four faces,
$b_1$, $b_2$, $b_3$ and $b_4$ (see
section~\ref{sec:combinatorics}). The local tiling near the ridges
$b_1\cap b_2$, $b_2\cap b_3$ and $b_3\cap b_4$ imply that the region
between the bisectors $\B(p_0,G_1^{\pm 1} p_0)$ is tiled by $E$,
$G_2G_1G_2(E)=G_2G_1(E)$ and $(G_2G_1G_2)^{-1}(E)=G_3^{-1}(E)$.

Note that none of the isometries mapping these three copies of $E$
fixes $p_1$, hence the only vertex cycle transformation for $p_1$ is
$G_1$, which is parabolic.

Now that we have checked cycles of ridges and boundary vertices, the
Poincar\'e polyhedron theorem shows that $E$ is a fundamental domain
for the action of $\Gamma$ modulo the action of $G_2$ (the latter
isometry generates the stabilizer of the center $p_0$ in
$\Gamma$). The main consequences will be drawn in
section~\ref{sec:presentation}.

We state the above result about cycles of boundary vertices in a
slightly stronger form.
\begin{prop}\label{prop:stab}
  The stabilizer of $p_1$ in $\Gamma$ is the cyclic group generated by
  $G_1$. The stabilizer of $q_1$ is generated by $G_2^{-1}G_3$.
\end{prop}

\subsection{Presentation}\label{sec:presentation}

The Poincar\'e polyhedron theorem (see section~\ref{sec:dirichlet})
gives the following presentation
$$
\langle G_1, G_2 | G_2^4, (G_1G_2)^3, (G_1G_2^2)^3\rangle
$$
or in other words, since $G_1G_2=G_2G_3$, 
$$
\langle G_2, G_3 | G_2^4, (G_2G_3)^3, (G_2G_3G_2)^3\rangle.
$$ 

It also gives precise information about the elliptic elements in the
group.
\begin{prop}\label{prop:nofixedpoint}
  Let $\gamma\in\Gamma$ be a non trivial torsion element. Then
  $\gamma$ has no fixed point in $\partial_\infty \ch2$.
\end{prop}
\begin{pf}
It follows from the Poincar\'e polyhedron theorem that any elliptic
element in $\Gamma$ must be conjugate to some power of a cycle
transformation of some cell in the skeleton of the fundamental
domain. This says that any elliptic element in the group must be
conjugate to a power of $G_2$ (which fixes the center of the Dirichlet
domain), a power of $G_1G_2$ (which preserves the ridge $b_1\cap b_2$)
or a power of $G_1G_2^2$ (which preserves the ridge $b_8\cap b_1$),
see section~\ref{sec:ridgeCycles}.

$G_1G_2$ and $G_1G_2^2$ are regular elliptic elements of order three,
so they do not fix any point in $\partial_\infty \ch2$ (nor do their
inverses). As for $G_2$, the only nontrivial, non regular elliptic
power is $G_2^2$, but this can easily be checked to be a reflection in
a point, so it is conjugate in ${\rm Bihol}(\mathbb{B}^2)$ to
$(z_1,z_2)\mapsto(-z_1,-z_2)$, which has no fixed point in the unit
sphere.
\end{pf}

\begin{prop} \label{prop:kernel}
  The kernel of $\rho_2$ is generated as a normal subgroup by $a^4$,
  $(at)^3$ and $(ata)^3$.
\end{prop}

\begin{pf}
  The fact that the three elements in the statement of the proposition
  are indeed in the kernel follows from the presentation and the fact that
  \begin{eqnarray*}
    & \rho_2(a)=G_2\\
    & \rho_2(b)=G_1^{-1}G_3\\
    & \rho_2(t)=G_3.
  \end{eqnarray*}

  We now consider the presentation
  $$
  \langle a,b,t | tat^{-1}=aba, tbt^{-1}=ab, a^4, (at)^3, (ata)^3\rangle.
  $$
  One can easily get rid of the generator $b$, since
  $$
  b = a^{-1}tat^{-1}a^{-1},
  $$ 
  and the other relation involving $b$ then follows from the other
  three relations. Indeed, one easily sees that $(at)^3=(ata)^3=1$ implies
  $(tat)^3=1$, and then
  $$ 
  t(\mathbf{a^{-1}tat^{-1}a^{-1}})t^{-1} = ta^{-1}ta(tat)^2 = ta^{-1}
  (ta)^2 t^2at = ta^{-1}\cdot a^{-1}t^{-1}\cdot t^2at$$ $$= t a^2 tat = ta
  t^{-1} a^{-1} = a (\mathbf{a^{-1}tat^{-1}a^{-1}}).
  $$

  In other words, the quotient group is precisely
  $$
  \langle a,t | a^4, (at)^3, (ata)^3\rangle,
  $$
  which is the same as the image of $\rho_2$.
\end{pf}

\section{Combinatorics at infinity of the Dirichlet domain} \label{sec:combinatoricsAtInfinity}

The next goal is to study the manifold at infinity, i.e. the quotient
of the domain of discontinuity under the action of the group. The idea
is to consider the intersection with $\partial \ch 2$ of a fundamental
domain for the action on $\ch 2$. Recall that we did not quite
construct a fundamental domain in $\ch 2$, but a fundamental domain
modulo the action of a cyclic group of order $4$ (generated by $G_2$).

\begin{figure}[htbp]
  \centering
  \epsfig{figure=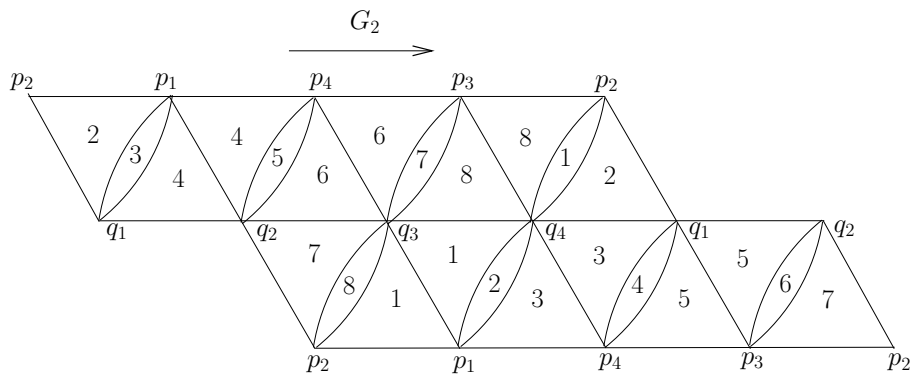, width=0.95\textwidth}
  \caption{The combinatorics of $\partial_\infty E$, which is a
    torus. We have split each quadrilateral components of the boundary
    faces into two triangles along an arc of $\C$-circle. Note that
    the polygons labelled~1 in this picture correspond to the 2-face
    illustrated in Figure~\ref{fig:b}.}
  \label{fig:torusFig8}
\end{figure}

\begin{figure}[htbp]
  \centering
 \epsfig{figure=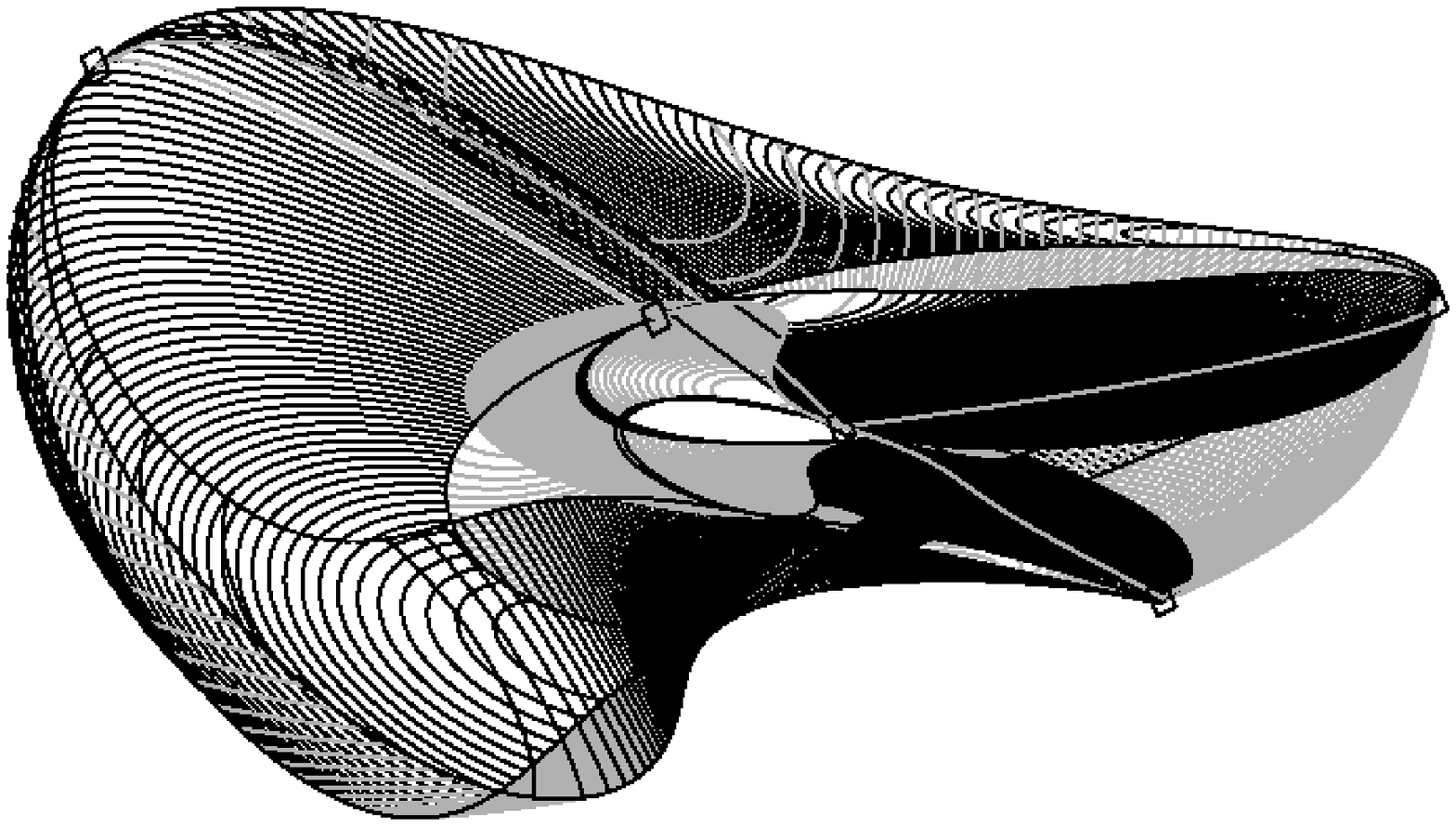, width = 0.8\textwidth}\\
  \epsfig{figure=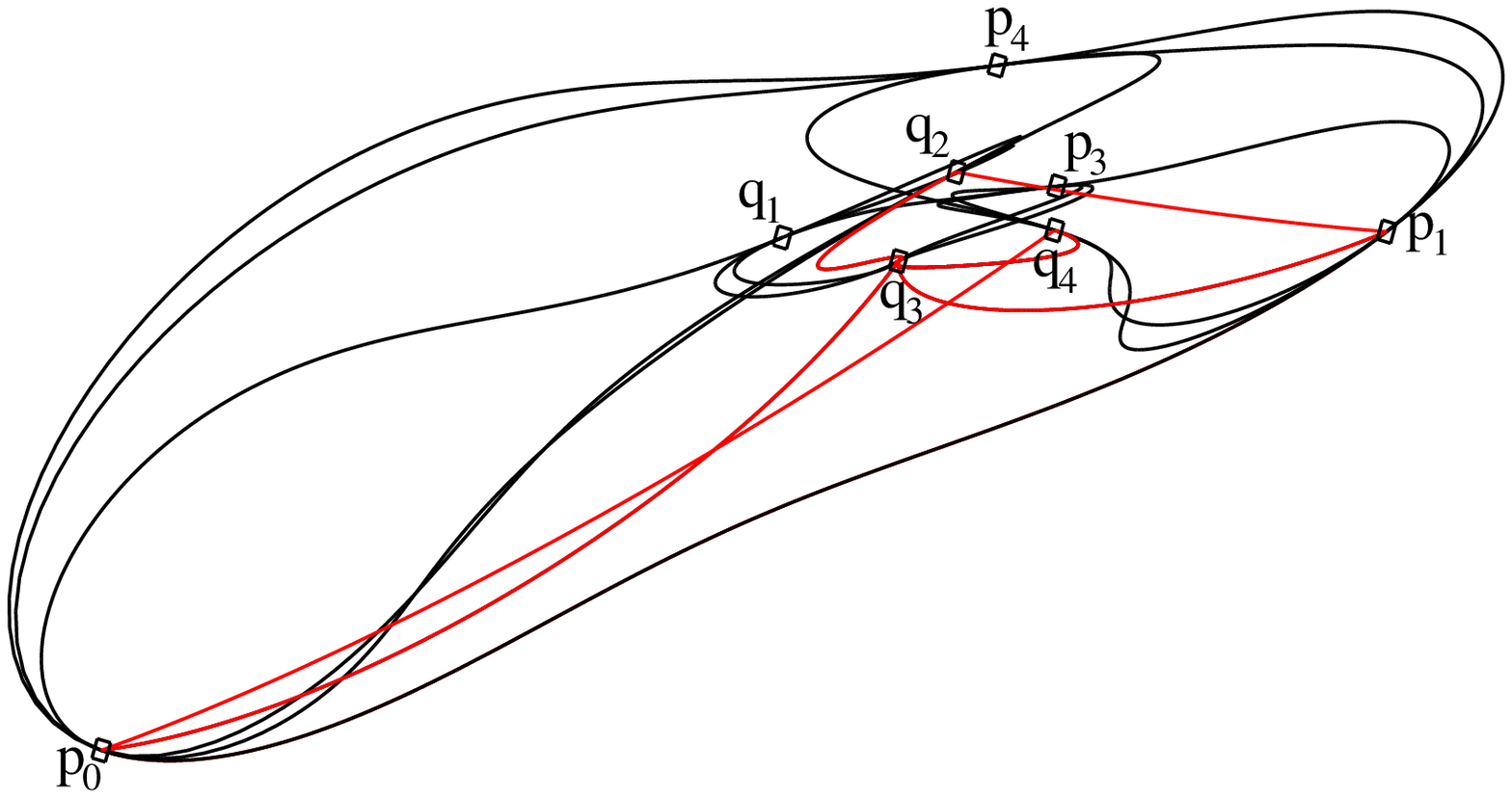, width = \textwidth}
  \caption{The solid torus $U$. On top, we have drawn all its
    $2$-faces, as well as its $1$-skeleton. On the bottom, only the
    $1$-skeleton with vertices labelled. These pictures are included
    only for motivational purposes, they are not needed in
    the proofs.}\label{fig:solidTorus}
\end{figure}

We start by describing the combinatorial structure of
$U=\partial_\infty E$, which is bounded by eight (pairwise isometric)
pieces of spinal spheres. A schematic picture of the boundary
$\partial U$ of $U$ in $\partial_\infty \ch 2$ is given in
Figure~\ref{fig:torusFig8}. The picture is obtained by putting
together the incidence information for each face, following the
results in section~\ref{sec:combinatorics}; we will use it as a
bookkeeping tool for the gluing of the eight faces. The
picture is by no means a realistic picture in complex hyperbolic space
(a more realistic view is given in Figure~\ref{fig:solidTorus}).
 
Note that it is clear from this picture that $\partial U$ is a torus,
and the fact that it is embedded in $\partial \ch 2$ follows from the
analysis of the combinatorics of $E$ given in the previous sections.
Figure~\ref{fig:solidTorus} makes it plausible that $U$ is a solid
torus. In fact a priori only one of the two connected components of
$S^3\setminus \partial U$ is a solid torus, the other may only be a
tubular neighborhood of a knot; in fact both sides are tori, because
one can produce two explicit simple closed curves with intersection
number one on $\partial U$, both trivial in $S^3$.  An alternative
argument for the fact that $U$ is a solid torus will be given below
(see Corollary~\ref{cor:solidtorus}).

\begin{rmk}
From the fact that $U$ is a solid torus, one can give a more direct
proof of the fact that the manifold at infinity of $\Gamma$ is the
figure eight knot complement. Indeed, Figure~\ref{fig:torusFig8} then
exhibits $U$ (with identifications on $\partial U$) as a 4-fold
covering of the figure eight knot. Rather than using this 4-fold cover
argument, we will divide $U$ into four explicit isometric regions, and
try modify the corresponding cell decomposition so that it is
combinatorially the same as the standard triangulation of the figure
eight knot complement. 
\end{rmk}

The next goal in our construction is to produce an explicit essential
disk in $U$ whose boundary is the curve on the left and right side of
Figure~\ref{fig:torusFig8}.
Note that $U$ is $G_2$-invariant simply because $E$ is so; the action
of $G_2$ on $\partial U$ is suggested on Figure~\ref{fig:torusFig8} by
the horizontal arrow. The rough idea is to use a fundamental domain
for the action of $G_2$ on $U$; the desired meridian would then be
obtained as one of the boundary components of this fundamental domain.

The Dirichlet domain has an arc in the boundary of a Giraud disk
between $q_1$ and $q_2$, which is in the intersection of the faces
$b_4$ and $b_5$. By Giraud's theorem (see~\cite{goldman}, p.~264),
there are precisely three bisectors containing that Giraud disk,
namely $\B_4$, $\B_5$, as well as
$$
\mathcal{C}=\B(G_1^{-1}p_0,G_2^{-1}G_3p_0).
$$ 
One way to get a fundamental domain for the action of $G_2$ on $U$
is to intersect $U$ with the appropriate region between $\mathcal{C}$
and $G_2\mathcal{C}$, namely
\small
$$
D=\left\{ z\in \C^3 : |\langle z,G_1^{-1}p_0\rangle |\leqslant |\langle
z,G_2^{-1}G_3p_0\rangle|, |\langle z,G_3p_0\rangle
|\leqslant |\langle z,G_2G_1^{-1}p_0\rangle|\right\}
$$ 
\normalsize
This turns out to give a slightly complicated fundamental domain
(in particular it is not connected).  We will only use $\mathcal{C}$
as a guide in order to get a simpler fundamental domain.

By construction, $\mathcal{C}$ contains $q_1$ and $q_2$. One easily
checks by direct computation that it also contains $p_2$, which is
given in homogeneous coordinates by $(0,0,1)$. To that end, one
computes
  $$
  \langle p_2, G_1^{-1} p_0\rangle=\langle p_2, G_2^{-1}G_3
  p_0\rangle=\frac{9+i\sqrt{7}}{4}.
  $$

\begin{figure}
    \centering
    \includegraphics[width=0.8\textwidth]{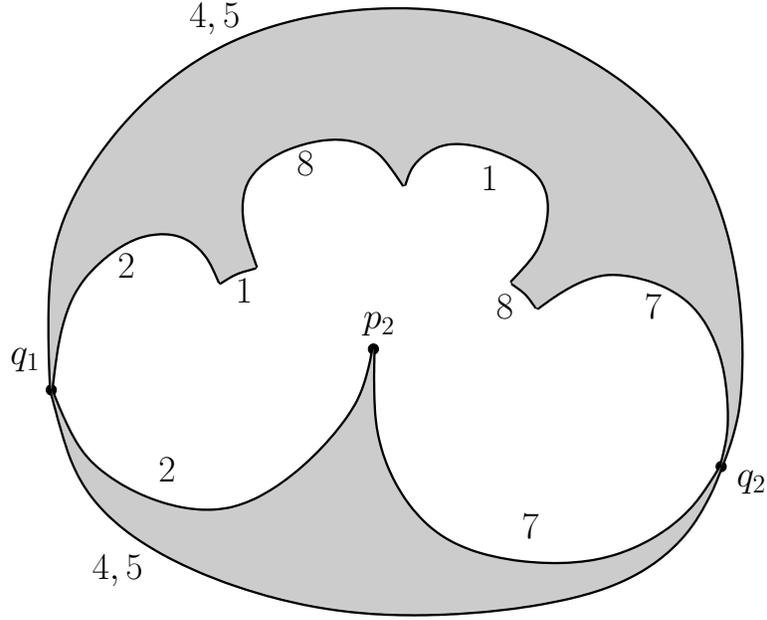}
    \caption{Combinatorics of the intersection of the spinal sphere
      $\partial_\infty \mathcal{C}$ with the solid torus $U$. The
      interior of this intersection has two components, one is a
      topological triangle with vertices $p_2$, $q_1$ and
      $q_2$.}\label{fig:picdisk}
  \end{figure}

One then studies the intersection of $\mathcal{C}$ with each face of
$U$ by using the techniques of section~\ref{sec:bisectors}. The only
difficulty is that the relevant bisectors are not all coequidistant
but their intersections turn out to be disks (this will be proved in
section~\ref{sec:noncoeq} of the appendix). The combinatorics of
$\mathcal{C}\cap E$ is illustrated in Figure~\ref{fig:picdisk}.
The picture suggests a natural way to choose an explicit parametrized
triangle $T$, with vertices $p_2$, $q_1$ and $q_2$ (and sides on the
appropriate bisector intersections, as indicated by labels in
Figure~\ref{fig:picdisk}).

\begin{figure}
  \centering 
  \epsfig{figure=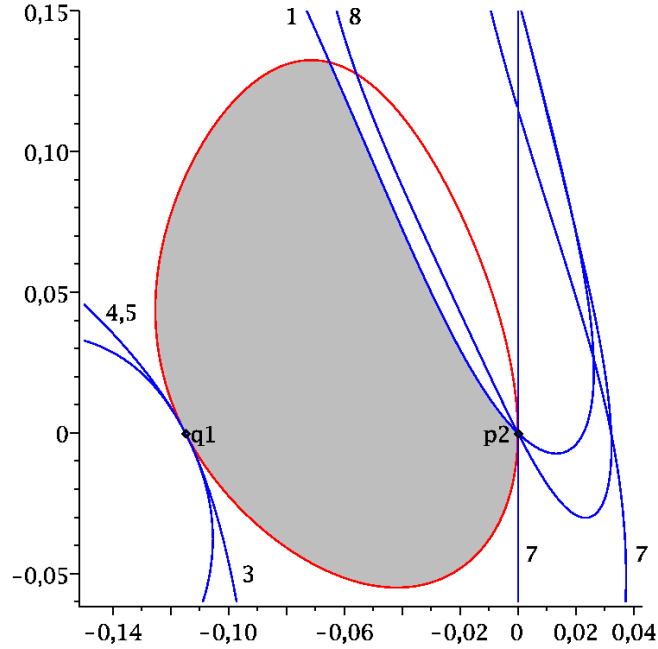,width=0.7\textwidth}
  \caption{Combinatorics of the intersection of $\mathcal{C}\cap\B_2 \cap E$.}
  \label{fig:2faces}
\end{figure}

Propositions~\ref{prop:deftau} and~\ref{prop:triangle} give a precise
definition of $T$ (their proof is quite computational, so we will give
them in the appendix, sections~\ref{sec:noncoeq}-\ref{sec:proofTriangle}).
\begin{prop}\label{prop:deftau}
  \begin{enumerate}
  \item $\partial_\infty(\B_4\cap\B_5)$ is a topological circle
    containing $q_1$, $q_2$ and $p_4$. We denote by $\tau_0$ the arc
    from $q_1$ to $q_2$ not going through $p_4$.
  \item $\partial_\infty(C\cap \B_7)$ is a topological circle
    containing $q_2$ and $p_2$, and only one of the two arcs of that
    circle from $q_2$ to $p_2$ is entirely contained in $U$; we write
    $\tau_1$ for that arc.
  \item $\partial_\infty(C\cap \B_2)$ is a topological circle
    containing $q_1$ and $p_2$, and only one of the two arcs of that
    circle from $p_2$ to $q_1$ is entirely contained in $U$; we write
    $\tau_2$ for that arc.
  \item The curve $\tau$ obtained by concatenating the arcs $\tau_0$,
    $\tau_1$ then $\tau_2$ from items (1), (2) and (3) is an embedded
    topological circle in $\partial_\infty C$.
  \end{enumerate}
\end{prop}
Item~(1) is obvious, since $\B_4\cap\B_5$ is a Giraud disk, and we
know which vertices lie on it (see
section~\ref{sec:combinatorics}). Items~(2) and~(3) follow from each
other by symmetry, we will only justify~(3).  The latter is made
plausible by Figure~\ref{fig:2faces}, which can be obtained using the
parametrizations explained in section~\ref{sec:bisectors}.

\begin{prop}\label{prop:triangle}
  The curve $\tau$ defined in Proposition~\ref{prop:deftau} bounds a
  unique triangle $T$ in $\partial_\infty C$ that is properly embedded
  in $U$. Moreover, $T$ and $G_2^2T$ are disjoint.
\end{prop}

An important consequence of Proposition~\ref{prop:triangle} is the
following.
\begin{cor}\label{cor:solidtorus}
  $U$ is an embedded solid torus in $\partial_\infty \ch 2$.
\end{cor}

\begin{pf}
The triangles $T$ and $G_2^2(T)$ split $U$ into two balls (they are
indeed balls because they are bounded by topological embedded
$2$-spheres), glued along two disjoint disks. From this it follows
that $U$ is a solid torus.
\end{pf}

In order to get a simple fundamental domain, we will modify the
meridian of Proposition~\ref{prop:triangle} slightly.
\begin{prop} \label{prop:homotopy}
  The side $\tau_1$ (resp. $\tau_2$) of $T$ is isotopic in the
  boundary to the arc of $\C$-circle joining these two points on the
  face $b_7$ (resp. $b_2$). Moreover, this isotopy can be performed
  so that the corresponding sides of the triangle $G_2(T)$ intersect
  the boundary of $T$ precisely in $q_2$.
\end{prop}

\begin{pf}
  The combinatorics of the face $b_2$ are combinatorially the same as
  Figure~\ref{fig:b}, but the pinch points are $p_1$ and $q_4$, and
  the other two vertices are $p_2$ and $q_1$ (see
  Figure~\ref{fig:torusFig8}). Since $\tau_2$ is contained in the face
  $b_2$ and contains no other vertex than $q_1$ and $p_2$, it remains
  in the interior of the quadrilateral component of $b_2$. In that
  disk component, any two paths from $q_1$ to $p_2$ are isotopic,
  hence all of them are isotopic to the path that follows the
  (appropriate) arc of the $\C$-circle between these two points.

  The argument for $\tau_1$ is similar. The fact that the isotopies
  for $T$ and $G_2(T)$ are compatible (in the sense that one can keep
  their sides disjoint throughout the isotopy) is obvious from the
  description of the combinatorics of $\partial U$, see
  Figure~\ref{fig:homotopy}.
  \begin{figure}
    \centering
    \epsfig{figure=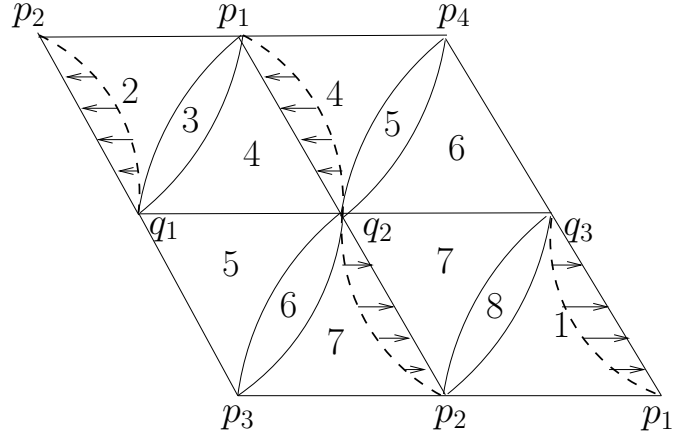, width=0.7\textwidth}
    \caption{Isotopy of part of the boundary of $T$ and $G_2(T)$
      towards an arc of a $\C$-circle.}\label{fig:homotopy}
  \end{figure}
\end{pf}

The upshot of the above discussion is that we have a convenient choice
of a meridian for the solid torus $U$, given by the concatenation of
the following three arcs
\begin{itemize}
\item The arc of $\C$-circle from $p_2$ to $q_1$ which is the
  boundary of a slice of the face $b_2$ (only one such arc is contained in
  the Dirichlet domain);
\item The arc of the boundary of the Giraud disk given by the
  intersection of the two bisectors $\B_4$ and $\B_5$, from $q_1$ to
  $q_2$ (there are two arcs on the boundary of this Giraud disk, we
  choose the one that does not contain $p_4$);
\item The arc of $\C$-circle from $q_2$ to $p_2$ which is the boundary
  of a slice of the bisector $\B_7$ (only one such arc is contained in
  the Dirichlet domain).
\end{itemize}
We denote this curve by $\sigma$.

\begin{prop}
  The curve $\sigma$ bounds a topological triangle $\tilde{T}$ which
  is properly contained in $U$. This triangle can be chosen so that
  $\tilde{T}\cap G_2 \tilde{T}$ consists of a single point, namely
  $q_2$.
\end{prop}

\begin{pf}
  This follows from the properties of $T$ and the isotopy of
  Proposition~\ref{prop:homotopy}.
\end{pf}

\section{The manifold at infinity}

The results from section~\ref{sec:combinatoricsAtInfinity} give a
simple fundamental domain for the action of $\Gamma$ in the domain of
discontinuity. For ease of notation, we denote $\tilde{T}$ simply by
$T$, see section~\ref{sec:combinatoricsAtInfinity} for how to obtain
this modified meridian for the solid torus $U$; recall that $U$ is by
definition the boundary at infinity $\partial_\infty E$ of the
Dirichlet domain $E$.

\begin{dfn}
  Let $D$ be obtained from the portion of $U$ that is between $T$ and
  $G_2(T)$.
\end{dfn}
By construction, this region has ten faces, eight coming from the
faces of the Dirichlet domain, and two given by $T$ and $G_2(T)$. For
each $k=1,\dots,8$, we denote by
$$
f_k=D\cap b_k
$$ 
the portion of $b_k$ that is inside $D$.

\begin{figure}
  \centering
  \epsfig{figure=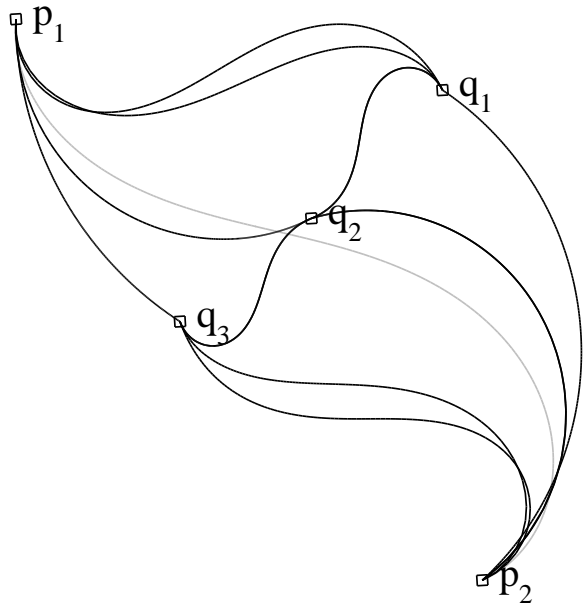, width=0.6\textwidth}
  \caption{A Heisenberg view of the $1$-skeleton of the fundamental
    domain $D$.}\label{fig:tetraNice}
\end{figure}

\begin{figure}
  \centering
  \epsfig{figure=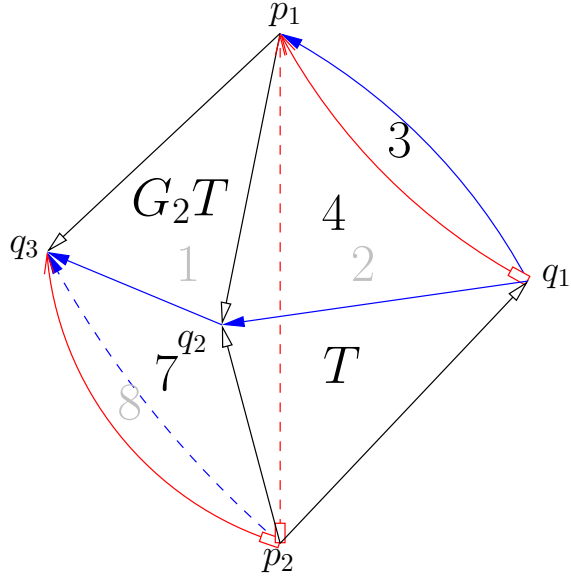, width=0.6\textwidth}
  \caption{The quotient manifold is homeomorphic to a ball with
    identifications on the boundary (one glues pairs of faces with
    matching arrows).}\label{fig:tetraCombinatorial}
\end{figure}

By construction, $D\cup G_2 D\cup G_2^2D\cup G_2^{-1}D$ is equal to
the solid torus $U=\partial_\infty E$. Since we have proved that $E$
tiles $\ch 2$, $U$ tiles $\partial_\infty \ch 2$ (in the sense that
eiher $U$ and $\gamma U$ coincide or $U\cap\gamma U$ has empty
interior). A Heisenberg view of the $1$-skeleton of $D$ is illustrated
in Figure~\ref{fig:tetraNice}, and a more combinatorial one, which we
will use later, is given in Figure~\ref{fig:tetraCombinatorial}.

The pictures we get are not quite the same as
Figure~\ref{fig:tetrahedra} (which is the one that usually appears in
the literature on the figure eight knot), but they are obtained from
it by taking the mirror image.

Note however that both oriented manifolds given by the usual or the
opposite orientation of the figure eight knot complement admit a
uniformizable spherical CR structure. Indeed, one can precompose the
developing map by an orientation-reversing automorphism of the figure
eight knot (hence the holonomy gets precomposed by the corresponding
automorphism of the fundamental group), see section~\ref{sec:mcg} for
more details.

Setting $V=\{p_1,\dots,p_4,q_1,\dots,q_4\}$, we also have that
$U^0=U\setminus V$ tiles the set of discontinuity $\Omega$ (indeed, it
follows from the Poincar\'e polyhedron theorem that the only fixed
points of parabolic elements in the group are conjugate to either
$p_1$ or $q_1$, see section~\ref{sec:poincare}).  We analyze the
quotient of $\Omega$ using the side pairings, which are given either
by the action of $G_2$ or by the side pairings coming from the
Dirichlet domain.

There are four side pairings, given in Table~\ref{tab:pairings}, three
coming from the Dirichlet domain, and one given by $G_2$. 
\begin{prop}
  The maps $G_1$, $G_2$, $G_3$ and $G_3G_1$ give side pairings of the
  faces of $D$, and map the vertices according to
  Table~\ref{tab:pairings}.
\begin{table}[htbp]
\centering
\begin{tabular}{rcl}
  $f_4$ & $\stackrel{G_1}{\longrightarrow}$ & $f_1$\\
  $p_1,q_2,q_1$ &  & $p_1,q_3,p_2$
\end{tabular}\\[0.5cm]
\begin{tabular}{rcl}
  $f_2$ & $\stackrel{G_3}{\longrightarrow}$ & $f_7$\\
  $p_1,p_2,q_1$ &  & $q_3,p_2,q_2$
\end{tabular}\\[0.5cm]
\begin{tabular}{rcl}
  $f_3$ & $\stackrel{G_3G_1}{\longrightarrow}$ & $f_8$\\
  $p_1,q_1$ &  & $q_3,p_2$
\end{tabular}\\[0.5cm]
\begin{tabular}{rcl}
  $T$ & $\stackrel{G_2}{\longrightarrow}$ & $G_2T$\\
 $p_2,q_1,q_2$ &  & $p_1,q_2,q_3$
\end{tabular}
\caption{The four side pairings, with their action on vertices. We
  denote by $f_k$ the part of $\partial_\infty b_k$ that is contained
  in $D$.}\label{tab:pairings}
\end{table}
\end{prop}

\begin{pf}
  The claim about $G_2$ holds by construction (see also
  Proposition~\ref{prop:tangency}). The ones about the other side
  pairings come from the Dirichlet domain (where an element $\gamma$
  maps the face associated to $\gamma^{-1}$ to the face associated to
  $\gamma$), see section~\ref{sec:sidePairings}.

  The claims about $G_3G_1$ follow from the previous ones, since
  $$
  G_3G_1(p_1)=G_3(p_1)=q_3
  $$
  and
  $$
  G_3G_1(q_1)=G_3(p_2)=p_2.
  $$
\end{pf}

We give a simple cut and paste procedure that allows us to identify
the quotient as the figure eight knot complement, and this will
conclude the proof of Theorem~\ref{thm:main}.

\begin{figure}
  \centering
  \epsfig{figure=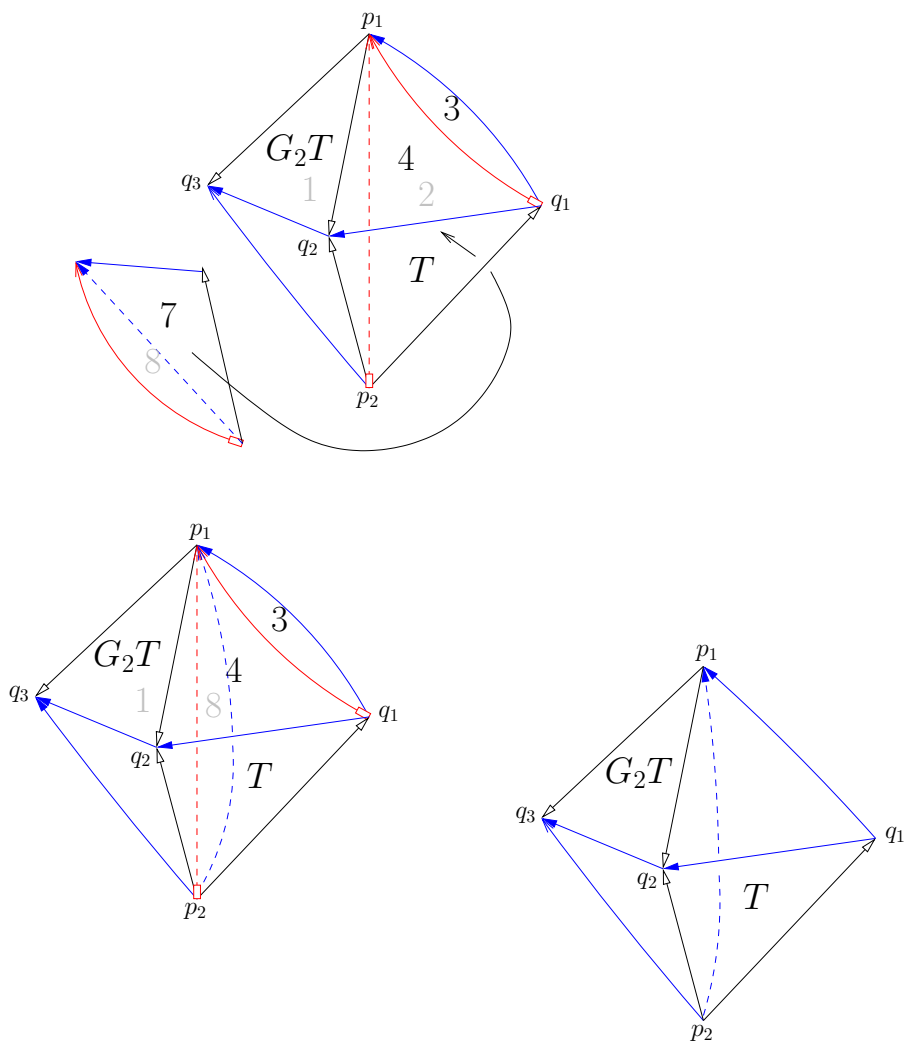, width=0.95\textwidth}
  \caption{Cut and paste instructions for recovering the usual
    two-tetrahedra decomposition of the figure eight knot
    complement.}\label{fig:cutAndPaste}
\end{figure}

The procedure is illustrated in Figure~\ref{fig:cutAndPaste}. We slice
off a ball bounded by $f_7$, $f_8$ as well as a triangle contained in
the interior of $D$, and move it in order to glue it to face $f_2$
according to the side pairing given by $G_3^{-1}$. Now we group faces
$f_1$ and $f_8$ on the one hand, and faces $f_4$ and $f_3$ on the
other hand, and observe that their side pairings agree to give the
identifications on the last domain in
Figure~\ref{fig:cutAndPaste}. This is the same as
Figure~\ref{fig:tetrahedra} (with the orientation reversed).

\section{Relationship between $\Gamma_2$ and $\Gamma_3$}\label{sec:thirdRep}

The goal of this section is to show that the groups $\Gamma_2$ and
$\Gamma_3$ are conjugate subgroups of $\pu21$.

We write
$$
G_1=\rho_2(g_1),\ G_2=\rho_2(g_2),\ G_3=\rho_2(g_3) 
$$
and
$$
A_1=\rho_3(g_1),\ A_2=\rho_3(g_2),\ A_3=\rho_3(g_3). 
$$

One can easily check that $A_1A_3^{-1}$ is regular elliptic element of
order $4$, hence it is tempting to take its isolated fixed point as
the center of a Dirichlet domain for $\Gamma_3$ (just like we did for
$\Gamma_2$, using the fixed point of $G_2$).

In fact it is easy to see that the corresponging Dirichlet domain is
isometric to that of $\Gamma_2$, and to deduce a presentation for
$\Gamma_3$, say in terms of the generators $M=A_1A_3^{-1}$ and
$N=A_1$:
$$
\langle M,N | M^4, (MN)^3, (MNM)^3\rangle
$$

With a little effort, these observations also produce an explicit
conjugacy relation between both groups. Denote by $P$ the following matrix:
$$
P=\left[\begin{matrix}
1 & 0 & 0 \\
\frac{-3-i\sqrt{7}}{4} & \frac{-5+i\sqrt{7}}{4} & 0\\
\frac{-1+i\sqrt{7}}{2} &  \frac{-1+i\sqrt{7}}{2} & 2
\end{matrix}\right]
$$

Then one easily checks (most comfortably with symbolic computation
software!) that
\begin{eqnarray*}
  &P^{-1} A_1 P = G_1^{-1}G_3G_1\\
  &P^{-1} A_3 P = G_3\\
\end{eqnarray*}
Note that the above two matrices generate $\Gamma_2$.
We will explain the precise relationship between the two representations
$\rho_2$ and $\rho_3$ in section~\ref{sec:mcg}.

\section{Action of $\textrm{Out}(\pi_1(M))$} \label{sec:mcg}

The main goal of this section is to explain the relationship between
the two representations $\rho_2$ and $\rho_3$, which turn out to
differ by precomposition with an outer automorphism of
$\pi_1(M)$. This is contained in the statement of
Proposition~\ref{prop:mcgaction}, where we analyze the action of the
whole outer automorphism group of $\pi_1(M)$.

\begin{figure}[htbp]
  \centering
  \epsfig{figure=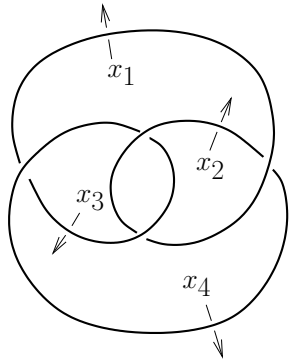, width=0.3\textwidth}
  \caption{A symmetric diagram for the figure eight knot - there are
    three planes of symmetry, one being the plane containing the
    projection.}\label{fig:symmetryKnot}
\end{figure}

We start by describing the outer automorphism group of $\pi_1(M)$ in
terms of explicit generators (it is well known that this group is a
dihedral group $D_4$ of order $8$). In fact $Out(\pi_1(M))$ can be
visualized purely topologically in a suitable projection of the figure
eight knot, for instance the one given in
Figure~\ref{fig:symmetryKnot}.

The Wirtinger presentation (see~\cite{rolfsen} for instance) is given
by
$$
\langle x_1,\dots,x_4\ |\ x_4x_1=x_3x_4,\ x_2x_3=x_3x_1,\ x_3x_2=x_2x_4,\ x_2x_1=x_1x_4\ \rangle.
$$
We eliminate $x_2$, then $x_3$ using
\begin{equation}\label{eq:x2x3}
x_2=x_1x_4x_1^{-1},\quad  x_3=x_2x_4x_2^{-1}
\end{equation}
and get
$$
\langle x_1,x_4\ |\ x_4[x_1^{-1},x_4]=[x_1^{-1},x_4]x_1\ \rangle.
$$ 
It will be useful to observe that with this presentation, we can express
$$
x_3=x_1x_4[x_1^{-1},x_4]x_1^{-1}=x_1[x_1^{-1},x_4]=x_4x_1x_4^{-1}.
$$

Of course the above presentation is the same as the one given in
section~\ref{sec:threereps} if we set
$$
x_1=g_3^{-1},\quad x_4=g_1^{-1}.
$$

\begin{figure}[htbp]
  \centering
  \epsfig{figure=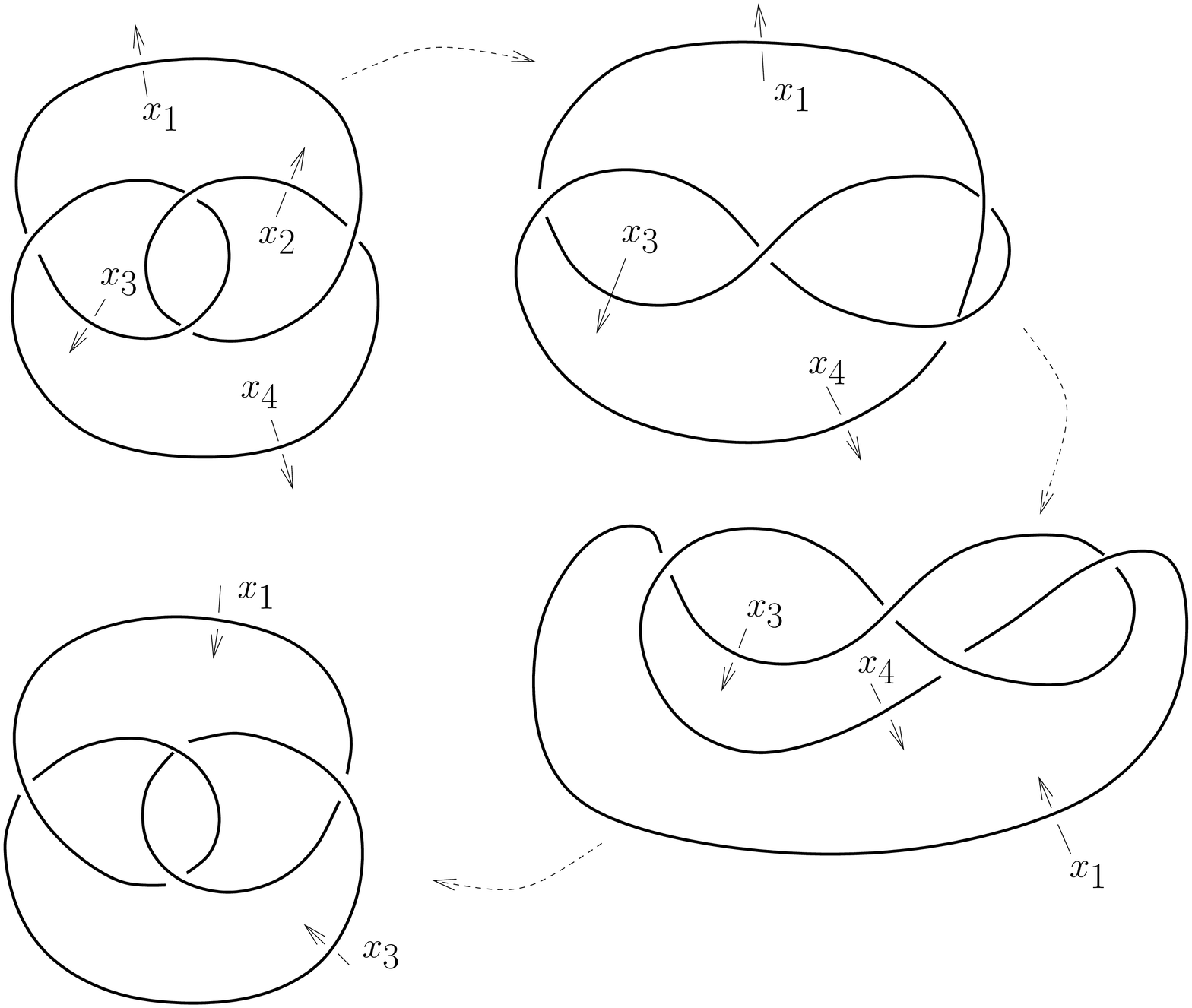, width=\textwidth}
  \caption{An isotopy from the figure eight knot to its mirror
    image.}\label{fig:isotopy}
\end{figure}

Using the Wirtinger presentation and an isotopy between the figure
eight knot and its mirror image, for instance as suggested in
Figure~\ref{fig:isotopy}, one can check that the automorphisms
described in Table~\ref{tab:mcggenerators} generate
$\textrm{Out}(\pi_1(M))$.

\begin{table}
\begin{eqnarray*}
  &\sigma:\left\{\begin{array}{l}
  g_1\mapsto g_3\\
  g_3\mapsto g_1
  \end{array}
  \right.\quad
  \iota:\left\{\begin{array}{l}
  g_1\mapsto g_1^{-1}\\
  g_3\mapsto g_3^{-1}
  \end{array}
  \right.
  \quad
  \tau:\left\{\begin{array}{l}
  g_1\mapsto g_1^{-1}g_3g_1\\
  g_3\mapsto g_3
  \end{array}
  \right.
\end{eqnarray*}
\caption{The three automorphisms $\sigma$, $\tau$, $\iota$ generate
  $Out(\pi_1(M))$.}\label{tab:mcggenerators}
\end{table}
Note that $\sigma$ and $\iota$ correspond to orientation-preserving
diffeomorphisms (and they generate a group of order $4$), whereas
$\tau$ reverses the orientation.

In what follows, for two representations $\rho$ and $\rho'$, we write
$\rho\sim\rho'$ when the two representations are conjugate. We start
with a very basic observation, valid for any unitary representation
(not necessarily with Lorentz signature).
\begin{prop}\label{prop:contragredient}
  Let $\rho:\pi_1(M)\rightarrow \un21$. Then $\rho\circ\iota\sim
  \overline{\varphi}^T$.
\end{prop}

\begin{pf}
  For any element $A$ of $\un21$,
  $$
  \overline{A}^TJA=J, 
  $$ 
  hence $A^{-1}=J^{-1}\ \overline{A}^TJ$ is conjugate to
  $\overline{A}^T$. 
\end{pf}

The precise relationship between $\rho_2$ and $\rho_3$ is as follows
(we only give the action of $\textrm{Out}(\pi_1(M))$ on $\rho_2$,
since the action on $\rho_3$ can easily be deduced from it).
\begin{prop}\label{prop:mcgaction}
  Let $\varphi\in\textrm{Out}(\pi_1(M))$. Then
  \begin{itemize}
    \item $\rho_2\circ\varphi\sim \rho_2$ if and only if $\varphi$ is
      trivial or $\varphi=\sigma\iota$.
    \item $\rho_2\circ\varphi\sim \overline{\rho}_2$ if and only if
      $\varphi=\sigma$ or $\iota$.
    \item $\rho_2\circ\varphi\sim \rho_3$ if and only if $\varphi=\tau$ or $\iota\tau$.
    \item $\rho_2\circ\varphi\sim \overline{\rho}_3$ if and only if
      $\varphi=\sigma\iota\tau$ or $\sigma\tau$.
  \end{itemize}
\end{prop}

\begin{pf}
  The fact that $\rho_2\circ \sigma\iota\sim \rho_2$ follows from the
  fact that $IG_1I=G_3^{-1}$, $IG_3I=G_1^{-1}$ (see
  section~\ref{sec:symmetry}). 

  One easily checks that
  \begin{equation}\label{eq:magic}
  G_1^T = G_3^{-1},\quad G_3^T = G_1^{-1}.
  \end{equation}
  Now the pair $G_1^{-1},G_3^{-1}$ is conjugate to
  $\overline{G}_1^T,\overline{G}_3^T$ (because the matrices preserve
  $J$), which is conjugate to
  $\overline{G}_3^{-1},\overline{G}_1^{-1}$ (by~(\ref{eq:magic})),
  which is conjugate to $\overline{G}_1,\overline{G}_3$ (by
  conjugation by $I$). This shows that $\rho_2\circ \iota\sim \overline{\rho}_2$.

  All that is left to prove is that $\rho_2\circ\tau\sim\rho_3$, and
  this was proved in section~\ref{sec:thirdRep}.
\end{pf}

\section{Appendix - sample calculations}\label{sec:appendix}

In this section we detail some of the computations that were mentioned
in previous sections of the paper (the general computational strategy,
and the geometric preliminaries are explained in
section~\ref{sec:bisectors}).
Throughout the appendix, we denote by $\hB_j$ denotes the extor in
projective space extending $\B_j$ (see~\cite{goldman} for a definition
and many properties of extors), and by $\overline{\B}_j$ the closure
of $\B_j$ in $\chb 2$. In other words,
$\overline{\B}_j=\B_j\cup\partial_\infty\B_j$. More generally,
$\widehat{O}$ denotes the extension to projective space of $O$, and
$\overline{O}$ denotes its closure in $\chb 2$.

\subsection{Pairs of bounding bisectors - proof of Lemma~\ref{lem:bisectorIntersections}} 
\label{sec:boundingpairs}

The center of the Dirichlet domain is given by \footnotesize
$$
p_0=(1,-\frac{3+i\sqrt{7}}{4},-1).
$$
\normalsize
Its relevant orbit points are given by
\footnotesize
$$
r_1=G_1 p_0=(\frac{3+i\sqrt{7}}{4}, \frac{1-i\sqrt{7}}{4},-1),\quad
r_2= G_3^{-1}p_0 = (1,\frac{1-i\sqrt{7}}{4},\frac{-3-i\sqrt{7}}{4})
$$
$$
r_3=G_2 r_1=(2,\frac{-1-i\sqrt{7}}{2},\frac{-3-i\sqrt{7}}{4}),\quad 
r_4= G_2 r_2 = (\frac{9-i\sqrt{7}}{4},\frac{-7-i\sqrt{7}}{4},-1)
$$
$$
r_5=G_2^2 r_1=(\frac{9-i\sqrt{7}}{4},\frac{-5-i\sqrt{7}}{2},-2),\quad
r_6= G_2^2 r_2 = (2,\frac{-5-i\sqrt{7}}{2},\frac{-9+i\sqrt{7}}{4})
$$
$$
r_7=G_2^{-1} r_1=(1,\frac{-7-i\sqrt{7}}{4},\frac{-9+i\sqrt{7}}{4}),\quad
r_8= G_2^{-1} r_2 = (\frac{3+i\sqrt{7}}{4},\frac{-1-i\sqrt{7}}{2},-2).
$$
\normalsize 

The Giraud torus $\OB_j\cap\OB_k$ can be parametrized by using the
techniques of section~\ref{sec:bisectors}. We start by proving
Lemma~\ref{lem:bisectorIntersections}.

In order to show that $\B_1\cap \B_2$ is a disk, it is enough to
exhibit a single point inside it, for instance
\begin{equation}\label{eq:x12}
X_{12}=(p_0-r_1)\boxtimes(p_0-r_2)=(\frac{1+i\sqrt{7}}{4}, \frac{3-i\sqrt{7}}{8}, -\frac{1+i\sqrt{7}}{4})
\end{equation}
does the job, since $\langle X_{12}, X_{12}\rangle=-3/4$.
Similarly, $\B_1\cap\B_3$ is a disk, because
$$
X_{13}=(p_0-r_1)\boxtimes(p_0-r_3)=(\frac{5+i\sqrt{7}}{8}, \frac{3-i\sqrt{7}}{8}, -\frac{1+i\sqrt{7}}{4})
$$
satisfies $\langle X_{13}, X_{13}\rangle=-1/2$.

In order to show that $\B_1\cap \B_5$ is empty, we parametrize the
Giraud torus $\OB_1\cap\OB_5$ by vectors of the form
$
( \overline{z}_1 p_0 - r_1 )\boxtimes (\overline{z}_2 p_0 - r_5), 
$ 
so that $V(z_1,z_2)$ is given by $v_0+z_1v_1+z_2v_2$, where
\footnotesize
$$
v_0 = ( \frac{9-3i\sqrt{7}}{8}, \frac{3+3i\sqrt{7}}{4}, 3 ),\quad 
v_1 = v_2 = (\frac{3-i\sqrt{7}}{8}, -\frac{1+i\sqrt{7}}{4}, -1),
$$ 
\normalsize 
see equation~\eqref{eq:paramTorus}. We then write out
$$
\langle V(z_1,z_2), V(z_1,z_2)\rangle=\Re(\mu(z_1)z_2)-\nu(z_1),
$$
where
$$
\mu(z_1)=\frac{5}{2}(\overline{z}_1-3),\quad \nu(z_1)=-\frac{55}{4}+\frac{15}{2}\Re(z_1).
$$
It is easy to verify that $\nu^2-|\mu|^2$ is always positive for
$|z_1|=1$, for instance by writing $z_1=x+iy$, and computing
$$
\nu(z_1)^2-|\mu|^2=225(2x-3)^2/16.
$$
Note that in order to get the previous formula, we have used the fact
that $|z_1|^2=x^2+y^2=1$.

\subsection{Proof of Lemma~\ref{lem:faceIntersections}} \label{sec:proof-faceint}

We first treat the proof of part~(1) of
Lemma~\ref{lem:faceIntersections}; even though, strictly speaking, it
will not be needed in the proof, we strongly suggest that the reader
keep Figure~\ref{fig:disks} in mind. We work only with $\B_1\cap
\B_7$, since $\B_1\cap \B_3$ can be deduced from it by symmetry.

The Giraud torus $\hB_1\cap \hB_7$ can be parametrized by vectors of
the form $(\overline{z}_1p_0-r_1)\boxtimes(\overline{z_2}p_0-r_7)$,
with $|z_1|=|z_2|=1$. In other words, we normalize it to be the
Clifford torus.

Explicitly, this can be written as $V(z_1,z_2)=v_0 +
z_1 v_1 + z_2 v_2$, where
\footnotesize
$$
v_0 = ( \frac{9-3i\sqrt{7}}{8}, \frac{-9+3i\sqrt{7}}{8}, \frac{15+3i\sqrt{7}}{8} )
$$
$$
v_1 = (-1, \frac {5+i\sqrt{7}}{4}, \frac{3-i\sqrt{7}}{8})
$$
$$
v_2 = ( \frac{-3+i\sqrt{7}}{8}, \frac{-1-i\sqrt{7}}{4}, -1)
$$ 
\normalsize
The Giraud disk inside the Clifford torus is described by imposing
that the above vector $V=V(z_1,z_2)$ be negative, i.e. $\langle
V,V\rangle<0$ which can be written as
\small
\begin{equation}\label{eq:bdy}
  \Re\left(7+\frac{-15-3i\sqrt{7}}{4}z_1+\frac{-15+3i\sqrt{7}}{4}z_2
         +\frac{i\sqrt{7}}{2}z_1\overline{z}_2\right)<0
\end{equation}
\normalsize
The equations of the intersection with $\hB_j$, $j=1,\dots,8$ are given by
\begin{equation}\label{eq:curve}
|\langle V(z_1,z_2),p_0 \rangle |^2 = |\langle V(z_1,z_2), r_j \rangle |^2, 
\end{equation}
and we write them in simplified form in Table~\ref{tab:equations} (by
simplified, we mean that we use $|z_1|=|z_2|=1$).
\begin{table}
\begin{tabular}{|c|l|}
\hline
$\B_1$ & $0$\\
$\B_2$ & $\Re( -43 + (-12+12i\sqrt{7}) z_1 + (33-3i\sqrt{7}) z_2 + (9-5i\sqrt{7}) z_1\overline{z}_2)\,/\,8$\\
$\B_3$ & $\Re( -{81} + {54} z_2 )\,/\,8$\\
$\B_4$ & $\Re( -{151} + (60+12i\sqrt{7}) z_1 + (60-12i\sqrt{7}) z_2 - (9+5i\sqrt{7}) z_1\overline{z}_2)\,/\,8$\\
$\B_5$ & $\Re( -{81} + {54} z_1 )\,/\,8$\\
$\B_6$ & $\Re( -{43} + (33+3i\sqrt{7}) z_1 - (12+12i\sqrt{7}) z_2 + (9-5i\sqrt{7}) z_1\overline{z}_2)\,/\,8$\\
$\B_7$ & $0$\\
$\B_8$ & $\Re( -16 + (15+3i\sqrt{7}) z_1 + (15-3i\sqrt{7})z_2 - (9+5i\sqrt{7}) z_1\overline{z}_2)\,/\,8$\\\hline
\end{tabular}
\caption{The equations of the intersection of each $\hB_j$ with the
  Giraud torus containing $\hB_1\cap \hB_7$ (the Giraud torus is given
  by the Clifford torus $|z_1|=|z_2|=1$).}\label{tab:equations}
\end{table}

\begin{prop}\label{prop:disjoint}
  For $j=3,4$ and $5$, $\B_j$ does not intersect $\B_1\cap\B_7$.
\end{prop}

\begin{pf}
  This was already proved for $\hB_3$ and $\hB_5$, since we proved
  Lemma~\ref{lem:bisectorIntersections} in
  section~\ref{sec:boundingpairs} (is says that $\B_1\cap\B_3$ and
  $\B_1\cap\B_5$ are empty). Alternatively, this can also be recovered
  from the equations given in Table~\ref{tab:equations}. For instance,
  the equation
  $$
  -81+27(z_2+\overline{z}_2)=0
  $$ 
  has only one solution given by $z_2=3/2$, which is not on the
  unit circle.

  We claim that the intersection with $\hB_4$ is empty as well. One way
  to see this is to write the equation in the form
  $$
  \Re(\mu(z_1) z_2) = \nu(z_1),
  $$ 
  which has a solution $z_2$ with $|z_2|=1$ if and only if
  $|\nu(z_1)|\leq |\mu(z_1)|$.

  In the case at hand, 
  \footnotesize
  $$
  \mu(z_1)=\frac{15-3i\sqrt{7}}{2} + \frac{-9+5i\sqrt{7}}{8}\overline{z}_1,\quad  \nu(z_1)=\frac{151}{8} - \Re\left(\frac{15+3i\sqrt{7}}{2} z_1\right).
  $$
  \normalsize
  One computes
  \footnotesize
  $$ 
  \nu^2-|\mu|^2=18193/64 - 2025 x/8 + 405\sqrt{7}y/8 - 42xy\sqrt{7}/2 + (209x^2 + 47y^2)/4 ,
  $$ 
  \normalsize
  where we have written $z_1=x+iy$. It is now standard 2-variable
  calculus to prove that this function is strictly positive on the unit
  disk.
\end{pf}

The extors $\hB_2$, $\hB_6$ and $\hB_8$ have 1-dimensional
intersection with the Giraud torus $\hB_1\cap\hB_7$. For the general
description of their (piecewise) parametrization by one spinal
coordinate, see section~\ref{sec:bisectors}. We explicit the
parametrization for $\hB_8$, since this will be needed in later
calculations.

\begin{figure}
  \centering
  \epsfig{figure=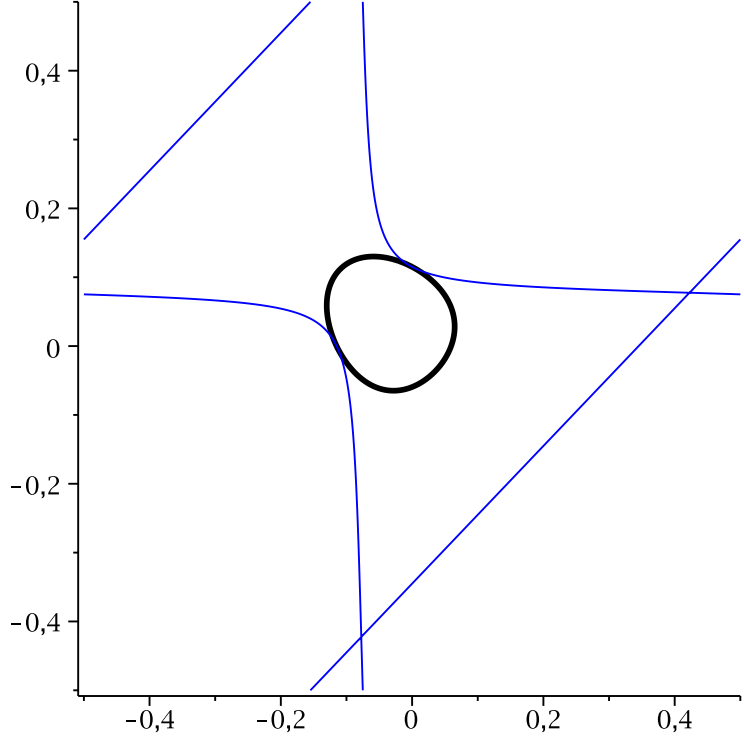, width=0.8\textwidth}
  \caption{The trace of $\B_8$ on the Giraud torus $\B_1\cap \B_7$, in
    terms of the log of the spinal coordinates (the bold oval is the
    boundary at infinity of the Giraud disk).}\label{fig:1_7_8}
\end{figure}

The equation for the trace on the Clifford torus of $\hB_8$ can be
written as
$$
  \Re(\mu(z_1) z_2)=\nu(z_1),
$$
where 
\small
$$
\mu(z_1) =  \frac{15-3i\sqrt{7}}{8}+ \frac{-9+5i\sqrt{7}}{8}\overline{z}_1,\quad \nu(z_1) =  2 - \Re \left( \frac{15+3i\sqrt{7}}{8}z_1 \right)
$$ 
\normalsize
Endpoints of the set of valid parameters are solutions of
$|\mu|^2=\nu^2$, which is a real expression involving $z_1$,
$\overline{z}_1$. Writing $z_1=x+iy$, we can write
$\nu^2-|\mu|^2=h(x,y)$ where
$$
h(x,y)= -1/2 - 45xy\sqrt{7}/32 - (31 x^2+193 y^2)/64.
$$ 
The endpoint of the parametrization are solutions of $h(x,y)=0$ that
satisfy $x^2+y^2=1$. The corresponding system has two solutions given
by $z_1=x+iy=\pm(5-i\sqrt{7})/4\sqrt{2}$, which have $\arg(z_1)/(2\pi)$
approximately equal to $-0.07745991$ and $0.42254009$ (compare with
the abscissas of the double points in
Figure~\ref{fig:1_7_8}).\label{endpoints}

Between these two values of the arguments, the sign of the
discriminant $\nu^2-|\mu|^2$ does not change, and it can easily be
checked that it is in fact nonpositive everywhere. In other words, the
formulae given in~\eqref{eq:param} parameterize the entire trace of
$\hB_8$ on the Clifford torus. The corresponding curve is depicted in
Figure~\ref{fig:1_7_8} (the figure is given only as a guide, it is not
needed in the proof).

Note that the curve seems to contain a straight line of slope
one. This is indeed the case, and it corresponds to a curve of the
form $z_2=\tau z_1$, for some complex number $\tau$ with
$|\tau|=1$. This straight line is actually contained in a complex
slice of the third bisector in Giraud's theorem, namely
$\B(r_1,r_7)$. Using the explicit form of the equation, plugging
$z_2=\tau z_1$, one finds a unique value of $\tau$ such that the
equation becomes trivial, namely
\begin{equation}\label{eq:tau}
  \tau=-\frac{9+5i\sqrt{7}}{16}.
\end{equation}
It is easy to see that this curve lies entirely outside complex
hyperbolic space. In fact substituting $z_2=\tau z_1$ in~\eqref{eq:bdy}
(and using $|z_1|=1$) yields a constant, namely $189/32$, which is
positive.

\begin{prop}\label{prop:tangent}
  For $j=2$, $6$ and $8$, $\B_j$ does not intersect
  $D=\B_1\cap\B_7$. In terms of their closures in $\chb 2$, we have
  the following:
  \begin{itemize}
    \item $\overline{\B}_2\cap \overline{D}=\{p_2\}$, which is the fixed point of $G_3$;
    \item $\overline{\B}_6\cap \overline{D}=\{q_3\}$, which is the fixed point of $G_1G_2^{-1}$;
    \item $\overline{\B}_8\cap \overline{D}=\{p_2,q_3\}$.
  \end{itemize}
  Moreover, (the extensions to projective space of) all these curves
  are tangent to $\partial_\infty D$ at every intersection point.
\end{prop}

\begin{pf}
  For $j=6$ and $7$, this follow from Proposition~\ref{prop:tangency}
  and Theorem~\ref{prop:unipotent} (since $\B_1$, $\B_6$,
  resp. $\B_2$, $\B_7$, have tangent spinal spheres).

  The statement about $j=8$ is a bit more difficult. We work in the
  Giraud torus normalized as the Clifford torus, which we write as
  $\widehat{D}$.
  We prove that the curves defined on $\widehat{D}$ by the equations
  for $\hB_2$ and $\hB_8$ are tangent at $p_2$ (a similar argument
  shows that the curves defined $\hB_6$ and $\hB_8$ are tangent at
  $q_3$).

  Recall that $p_2=(0,0,1)$, which we now need to write in the spinal
  coordinates $(z_1,z_2)$ for $\widehat{D}$. This is done by solving
  $\langle p_2, \overline{z}_1p_0-r_1\rangle=0$ for $z_1$, and
  $\langle p_2, \overline{z}_2p_0-r_2\rangle=0$ for $z_2$. Explicit
  calculation shows $p_2$ is given in spinal coordinates by
  $$
  (z_1,z_2)=(\frac{3-i\sqrt{7}}{4},1).
  $$  
  All equations in Table~\ref{tab:equations} have the form $f=0$ where
  $$
  f(z_1,z_2) = 2 \Re( a_0 + a_1 z_1 + a_2 z_2 + a_{12}z_1\overline{z}_2 ). 
  $$
  Since we are interested in the solution set only on the Clifford
  torus, we write $z_1=e^{it_1}$ and $z_2=e^{it_2}$ for real
  $t_j$. the gradient of $f$, seen as a function of $(t_1,t_2)$ is
  then given by
  $$
  \frac{\partial f}{\partial t_1} = - 2 \Im(a_2 z_2 + \overline{a}_{12} \overline{z}_1z_2),\qquad  
  \frac{\partial f}{\partial t_2} = - 2 \Im(a_2 z_2 + \overline{a}_{12} \overline{z}_1z_2).
  $$
  
  This gives
  \small
  $$
  \nabla f_2(\frac{3-i\sqrt{7}}{4},1) = (-\frac{3\sqrt{7}}{4},-\frac{3\sqrt{7}}{8}), \qquad
  \nabla f_8(\frac{3-i\sqrt{7}}{4},1) = (\frac{3\sqrt{7}}{8},\frac{3\sqrt{7}}{16}), 
  $$ 
  \normalsize
  where $f_j$ denotes the equation of $\B_j\cap \widehat{D}$, see
  Table~\ref{tab:equations}. This shows the needed tangency.

  It follows from Proposition~\ref{prop:unipotent} that $\hB_2\cap
  \widehat{D}$ is tangent to $\partial_\infty D$ at $p_2$. From the
  previous computation, we see that $\hB_2\cap \widehat{D}$ is also
  tangent to $\partial_\infty D$ at $p_2$.

  We now argue that $\hB_8\cap \widehat{D}=\{p_2,q_3\}$. Even though
  this is quite clear from the picture of the parametrized curve, we
  give a computational argument that does not rely on visual aids.

  We have explicit equations for $\partial_\infty D$ and $\hB_8$,
  namely~\eqref{eq:bdy} (with the inequality replaced by an equality)
  and~\eqref{eq:curve}. Writing out $z_j=x_j+iy_j$ for real $x_j,y_j$,
  the intersection $\partial_\infty D\cap \hB_8$ is described by the
  solutions of the system \small
  $$
  \left\{\begin{array}{l}
  15(x_1+x_2)+3\sqrt{7}(y_2-y_1)+2\sqrt{7}(x_2y_1-y_2x_1)=28\\
  15(x_1+x_2)+3\sqrt{7}(y_2-y_1)-9(x_1x_2+y_1y_2)+5\sqrt{7}(x_2y_1-y_2x_1)=16\\
  x_1^2+y_1^2=1\\
  x_2^2+y_2^2=1
  \end{array}
  \right.
  $$ 
  \normalsize 
  One checks that this has exactly two solutions, given
  by $(z_1,z_2)=(1,\frac{3+i\sqrt{7}}{4})$ (this corresponds to $q_3$)
  and $(z_1,z_2)=(\frac{3-i\sqrt{7}}{4},1)$ (this corresponds to
  $p_2$).

  Recall that $\widehat{D}\cap\widehat{\B}_8$ contains a diagonal
  component, given by $z_2=\tau z_1$ with $\tau$ as
  in~\eqref{eq:tau}. Recall that $\widehat{D}\cap\widehat{\B}_8$ has
  two double points, which were computed on
  page~\pageref{endpoints}. Away from these two endpoints, for a given
  $z_1\in S^1$, there is precisely one point $(z_1,z_2)$ in
  $\widehat{D}\cap\widehat{\B}_8$ that is \emph{not} in the diagonal
  component. The closure of that component (obtained by adding the two
  double points), gives an embedded topological circle in
  $\widehat{D}$. Since its only contact points with $\partial_\infty
  D$ are the two tangency points, we know this circle lies entirely
  outside $D$.
\end{pf}

This finishes the proof of part~(1) of
Lemma~\ref{lem:faceIntersections}. Part~(2) is very similar; by
symmetry, it is enough to consider $\B_1\cap\B_2$. 

As in the case of $\B_2\cap\B_7$, one finds all the intersections of
the Giraud torus $\OB_1\cap\OB_2$ with every $\OB_k$ ($k\neq 1,2$),
and checks that the only ones are given by $p_1$, $p_2$ and
$q_4$. This shows that $\B_1\cap\B_2\cap E$ is either empty or all of
$\B_1\cap\B_2$. One shows that it is a disk simply by finding one
point in it, for instance the point in~\eqref{eq:x12} is easily seen
to be inside $E$ by computing six inequalities.

This finishes the proof of Lemma~\ref{lem:faceIntersections}, hence
also of part~(1) of Proposition~\ref{prop:comb-b1}. Part~(2) will be
proved in section~\ref{sec:easy-spinal}.

\subsection{Spinal sphere of $\B_1$ - proof of Proposition~\ref{prop:comb-b1}(2)} \label{sec:easy-spinal}

In this section, we justify Proposition~\ref{prop:comb-b1}(2); in
other words, we justify the picture given in Figure~\ref{fig:b}.

We start by giving explicit coordinates on $\B_1=\B(p_0,G_1 p_0)$.  We
choose coordinates for $\ch 2$, seen as the unit ball $\mathbb{B}^2\subset\C^2$, where
the midpoint of the segment $[p_0,r_1]$ is taken to be at the origin of the ball
(as in section~\ref{sec:boundingpairs}, we write $r_1=G_1 p_0$). Since
$\langle p_0, r_1\rangle$ is real (and $\langle p_0,p_0\rangle=\langle
r_1,r_1\rangle$), the midpoint is given by $p_0+r_1$, and an
orthogonal vector spanning the complex spine is given by $p_0-r_1$.

We normalize these vectors to have unit norm, so we take
\footnotesize
$$
v_0=(p_0+r_1)/\sqrt{5}=(\frac{7+i\sqrt{7}}{4\sqrt{5}},\frac{-1-i\sqrt{7}}{2\sqrt{5}},-\frac{2}{\sqrt{5}}),
$$
$$
v_1=p_0-r_1=(\frac{1-i\sqrt{7}}{4},-1,0),
$$
$$
v_2=( \frac{-3+i\sqrt{7}}{4\sqrt{5}}, -\frac{1+i\sqrt{7}}{2\sqrt{5}}, -\frac{2}{\sqrt{5}}).
$$ 
\normalsize
The last vector is chosen so that $v_0,v_1,v_2$ is a standard Lorentz
basis, i.e. if $P$ denotes the corresponding base change matrix,
$$
P^* J P = \left(\begin{matrix}
  -1 & 0 & 0 \\
  0 & 1 & 0\\
  0 & 0 & 1
\end{matrix}\right).
$$ 
We now work in $\C^2$, with affine coordinates $u_1=z_1/z_0$,
$u_2=z_2/z_0$, where the $z_j$ denote coordinates in the basis
$v_0,v_1,v_2$; the complex hyperbolic plane $\ch 2$ is then simply
given by the unit ball $|u_1|^2+|u_2|^2<1$.

The ball coordinates for $p_0$ and $r_1$ are given by $(\pm
1/\sqrt{5},0)$, and the bisector $\B_1$ has a very simple equation,
namely 
$$
  \Re(u_1)=0,
$$
so the bisector can simply be thought of as the unit ball in $\R^3$,
when using coordinates $(t_1,t_2,t_3)$ for a point in $\mathbb{B}^2$ of the
form
$$
  (it_3,t_1+it_2).
$$
Here we have choosen the real spine of $\B_1$ to be given by the last
coordinate axis.

The equation of the intersection of a bisector $\B_j=\B(p_0,r_j)$ for
some $j>1$ is obtained simply by writing $r_j$ in the new
basis. In fact the equation has the form
\begin{equation}\label{eq:otherbisector}
  |\langle Z,P^{-1}p_0\rangle|^2 = |\langle Z,P^{-1}r_j\rangle|^2  
\end{equation}
where one takes $Z=(1,it_3,t_1+it_2)$.

We write $g_j$ for the equation of $\hB_2\cap \partial_\infty\B_1$ in
the coordinates $t_j$ for $\partial_\infty\B_1$ described
above. According to previous discussions (see
section~\ref{sec:combinatorics}), we only need to consider the
bisectors $\hB_2$ and $\hB_8$. The affine coordinates of $r_2$ and $r_8$
are given by
\small
$$
(0,\frac{-9+5i\sqrt{7}}{24}),\quad (0,\frac{2}{3}),
$$
\normalsize
respectively.

\begin{figure}
\centering
\epsfig{figure=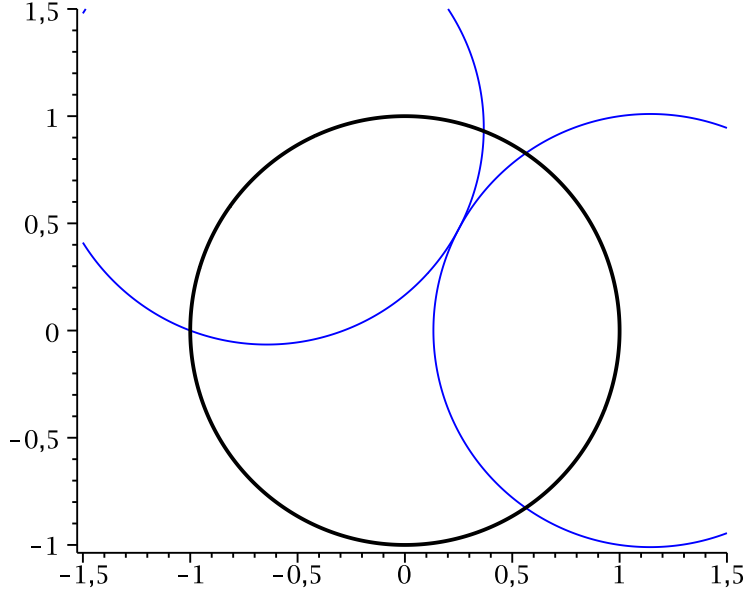,width=0.8\textwidth}
\caption{When $\B_1$ is normalized to be the unit ball with real spine
  given by the $t_3$-axis, $\B_1\cap \B_2$ and $\B_1\cap \B_8$ project
  to circles in the $(t_1,t_2)$-plane.}\label{fig:circles}
\end{figure}

We consider the intersection of $\hB_j$ with $\partial_\infty \B_1$,
the latter being given by the unit sphere
$t_1^2+t_2^2+t_3^2=1$. Computationally, we take the resultant $h_j$ of
$g_j$ and $t_1^2+t_2^2+t_3^2-1$ with respect to $t_3$.  For $j=2$ and
$8$ we get
\footnotesize
$$
h_2(t_1,t_2)=\frac{21}{20}\{(t_1+\frac{9}{14})^2+(t_2- \frac{5}{2\sqrt{7}})^2-\frac{50}{49}\},\quad
h_8(t_1,t_2)=\frac{21}{20}\{(t_1-\frac{8}{7})^2-t_2^2-\frac{50}{49}\}.
$$ 
\normalsize
The equations $h_2=0$ and $h_8=0$ define two cylinders in $\R^3$,
that project to a pair of tangent circles. The point of tangency of
the projections is given by $(1/4,5\sqrt{7}/28)$, as illustrated in
Figure~\ref{fig:circles}.
The inequalities defining the Dirichlet domain correspond to $g_j$
being negative. In particular, points in the interior of the Dirichlet
domain are the points in the unit ball that project outside both these
circles. 

It follows from the analysis in section~\ref{sec:combinatorics} and
the results in section~\ref{sec:proof-faceint} that $\partial_\infty
b_1$ is bounded only by the two curves corresponding to the
intersections with $\B_2$ and $\B_8$ (both of these curves are traces
on the $\partial_\infty\ch 2$ of Giraud disks). This finishes the
proof of Proposition~\ref{prop:comb-b1}.

\subsection{The intersection $C\cap\B_2$ is a disk} \label{sec:noncoeq}

In this section we consider $C\cap \B_2$, where
$C=\B(G_1^{-1}p_0,G_2^{-1}G_3 p_0)=\B(r_4,r_5)$. We will show that it
is a disk.

Note that these the bisectors $C$ and $\B_2$ do not share any complex
slice, i.e. their extended real spines do not intersect. This amounts
to saying that the circles $\overline{z}_1r_4-r_5$, $|z_1|=1$ and
$\overline{z}_1p_0-r_2$, $|z_2|=1$ do not intersect. 

One way to see this is to compute the intersection of their extended
complex spine, which can be represented by 
\small
$$
s=(-3\sqrt{7}-5i,\ 4\sqrt{7}+10i,\ 4\sqrt{7}),
$$ 
\normalsize
and to note that this vector satisfies $\langle s,s \rangle=44>0$. This point is
on the real spine of $C$ if and only if there exists a $z_1\in S^1$
such that $\langle s, \overline{z}_1r_4-r_5\rangle = 0$. The latter
can only happen if $z_1=(9+15i\sqrt{7})/46$, but this does not have
modulus one. Similarly one checks that $s$ is not on the real spine of
$\B_2$.

Now the intersection $C\cap \B_2$ can be parametrized by vectors of
the form $(\overline{z}_1 r_4-r_5)\boxtimes(\overline{z}_2 p_0 -
r_2)$. Such vectors have negative norm if and only if
\begin{equation}\label{eq:munu2}
  \Re(\mu(z_1) z_2) < \nu(z_1)
\end{equation}
where 
\footnotesize
$$
\mu(z_1)= ( -39 - 3i\sqrt{7} +(9+3i\sqrt{7})z_1 + 18\overline{z}_1 ) / 4,\quad
\nu(z_1)= - 15 + \Re\{ (24+3i\sqrt{7})z_1 \}/2.
$$
\normalsize
In order to analyze the number of connected components of the
intersection, we search for values of $z_1$ where the discriminant
vanishes. Writing $z_1=x+iy$, the discriminant $\nu(z_1)^2-|\mu(z_1)|^2$ becomes
\footnotesize
$$
\delta(x,y)=(1413 - 1764x + 216\sqrt{7}y +351(x^2-y^2) - 180\sqrt{7}xy)/8,
$$ 
\normalsize The system $\delta(x,y)=x^2+y^2-1=0$ has exactly two
solutions, one given by $(x,y)=(1,0)$, and the other one given by the
single real root of each of the polynomials $2221x^3-7103x^2+7411x-2473$,
$392\sqrt{7} +2268y + 1024\sqrt{7}y^2+2221y^2$. An approximate value
of $(x,y)$ is $(0.70552301,-0.70868701)$.

In fact only the number of solutions interests us; $z_1=e^{2\pi it_1}$
give nontrivial intervals of values of $z_2$ when $t_1^{min}<t_1<0$,
where $t_1^{min}=-0.12535607\dots$. For each such $z_1$, there
is only an interval of values of $z_2$ satisfying~\eqref{eq:munu2}, hence
$C\cap \B_2$ is a disk.

\subsection{Proof of Proposition~\ref{prop:deftau}(3)} \label{sec:proofdeftau3}

We consider the segment $\tau_2$, which corresponds to the bottom
segment from $q_1$ to $p_2$ shown on Figure~\ref{fig:2faces}. We prove
that it is contained in the (boundary at infinity of the) Dirichlet
domain; this will prove Proposition~\ref{prop:deftau}, since one
easily shows that the top arc of Figure~\ref{fig:2faces} is not
contained in $U$, simply by picking one point just above $p_2$.

It is enough to find all intersection points of $C\cap \B_2\cap\B_j$
for $j\neq 2$, and to show that none of them is in (the interior of)
the bottom segment; note that, in our coordinates, the bottom segment
is characterized by the fact that $\arg(z_2)<0$.

The (finite) list of points in $C\cap \B_2\cap\B_j$ can be obtained by
using Groebner bases. For instance, for $j=1$, the intersection points are
given by the solutions of the system \footnotesize
$$
\left\{\begin{array}{l}
-\frac{99}{8} -\frac{3\sqrt{7}}{2} y_1 +\frac{21}{2} x_1 -\frac{3}{2}x_1
   + \frac{15\sqrt{7}}{8}(x_1y_2 +y_1x_2) + \frac{27}{8}(x_1x_2-y_1y_2)=0\\
15 - 12 x_1 + \frac{3\sqrt{7}}{2} y_1 - \frac{39}{4} x_2 +
  \frac{3\sqrt{7}}{4}y_2 + \frac{27}{4}x_1x_2 + \frac{9}{4}y_1y_2
  -\frac{3\sqrt{7}}{4} (x_1y_2+x_2y_1)=0\\
x_1^2+y_1^2=1\\
x_1^2+y_2^2=1
\end{array}
\right.
$$ 
\normalsize
where we have split $z_1=x_1+iy_1$ and $z_2=x_2+iy_2$ into their real
and imaginary parts. This system has precisely two solutions, one
given by $(z_1,z_2)=(1,0)$, and the other one with
$$
\arg(z_1)/(2\pi)\approx  -0.06508170,\quad \arg(z_2)/(2\pi)\approx 0.13166662
$$

For $j=3$, the result follows from direct calculations in a similar
vein (using Groebner bases in order to solve the system). The
intersection of $\B_3$ is tangent to $\partial_\infty (C\cap\B_2)$, so
one gets a single intersection point, corresponding to $q_1$.

For $j=5$ or $7$, no computation is needed; we already know
that $\B_2\cap\B_5=\{q_1\}$ and $\B_2\cap \B_7=\{p_2\}$, since the
corresponding bisectors have tangent spinal spheres (see
section~\ref{sec:statement}).

\begin{rmk}
The intersections $C\cap \B_2\cap\B_j$ can also be handled by using
coequidistant pairs of bisectors, by writing the equation of the trace
of $C$ on $\B_2\cap\B_j$. 
\end{rmk}

\subsection{Proof of Proposition~\ref{prop:deftau}(4)} \label{sec:proofdeftau4}

In this section, we prove that the curve $\tau$ from
Proposition~\ref{prop:deftau} is an embedded topological circle in
$\partial_\infty C$. We also give explicit parametrizations of its
sides $\tau_0$, $\tau_1$, $\tau_2$, which are used to draw the
pictures in section~\ref{sec:proofTriangle}.

We start by parametrizing $\partial_\infty C$; we choose coordinates
for $\ch 2$ (seen as the unit ball $\mathbb{B}^2$) where the midpoint
of $[r_4,r_5]$ is at the origin (such a normalization was already
discussed in section~\ref{sec:easy-spinal}). A possible base change
matrix is given by
  \begin{equation}\label{eq:basechange2}
  P = \left(\begin{matrix}
    \frac{9-i\sqrt{7}}{2\sqrt{5}}     &          0            & \frac{-17+3i\sqrt{7}}{4\sqrt{5}} \\
    \frac{-17-3i\sqrt{7}}{4\sqrt{5}}  & \frac{3+i\sqrt{7}}{4}  & \frac{9+i\sqrt{7}}{2\sqrt{5}}\\
    -\frac{3}{\sqrt{5}}              &          1            &      \frac{2}{\sqrt{5}}
 \end{matrix}\right).
  \end{equation}
As in section~\ref{sec:easy-spinal}, we parametrize the spinal sphere
$\partial_\infty C$ as the unit sphere in $\R^3$ with coordinates
$t_j\in \R$, where $(t_1,t_2,t_3)$ corresponds to
  $
  (1,it_3,t_1+it_2).
  $ 
In these coordinates, the equations for the intersection of the
$\hB_j$ with $\partial_\infty C$ are then computed explicitly to be
those in Table~\ref{tab:equations2} (we obtain them by
simplifying~\eqref{eq:otherbisector}, using $t_1^2+t_2^2+t_3^2=1$).

The vertices of the triangle $T$ are given in
Table~\ref{tab:tvertices}.
  \begin{table}
  \tiny
  \begin{tabular}{|c|c|}
    \hline
    $\B_1$&$-69/10-6t_2t_3/\sqrt{5} + 66t_1/5 - 123(t_1^2+t_2^2)/20$\\    
    $\B_2$&$-24/5-39t_2t_3/8\sqrt{5}+261t_1/40 + 3\sqrt{7} t_2/8 - 3 t_3\sqrt{7}/\sqrt{5}-9(t_1^2+t_2^2)/5
            +21 t_1t_3\sqrt{7}/8\sqrt{5}$\\
    $\B_3$&$-21/10-33t_2t_3/8\sqrt{5} + 27 t_1/40 - 3 \sqrt{7}t_2/8 - 3 t_3\sqrt{7}/2\sqrt{5}+9(t_1^2+t_2^2)/10
            +3 t_1t_3\sqrt{7}/8\sqrt{5}$\\
    $\B_4$&$ 3/10-12 t_1/5 + 21(t_1^2+t_2^2)/20$\\       
    $\B_5$&$ 3/10-12 t_1/5 + 21(t_1^2+t_2^2)/20$\\       
    $\B_6$&$-21/10+33t_2t_3/8\sqrt{5} + 27 t_1/40 - 3 \sqrt{7}t_2/8 + 3 t_3\sqrt{7}/2\sqrt{5}+9(t_1^2+t_2^2)/10
            -3 t_1t_3\sqrt{7}/8\sqrt{5}$\\
    $\B_7$&$-24/5+39t_2t_3/8\sqrt{5}+261t_1/40 + 3\sqrt{7} t_2/8 + 3 t_3\sqrt{7}/\sqrt{5}-9(t_1^2+t_2^2)/5
            -21 t_1t_3\sqrt{7}/8\sqrt{5}$\\
    $\B_8$&$ -69/10+6t_2t_3/\sqrt{5} + 66t_1/5 - 123(t_1^2+t_2^2)/20$\\\hline    
  \end{tabular}
  \normalsize
  \caption{The equations of $\hB_j$ in $\partial_\infty C$, for $j=1,\dots,8$.}\label{tab:equations2}
  \end{table}
  \begin{table}
  \begin{tabular}{c|c|c}
    Vertex & $(t_1,t_2,t_3)$                   & $j$ such that $f_j=0$\\\hline
    $p_2$ & $(\frac{87}{88},\frac{5\sqrt{7}}{88},0)$             & 1,2,7,8\\
    $q_1$ & $(\frac{1}{4},\frac{5\sqrt{7}}{28},-\sqrt{\frac{5}{7}})$   & 2,3,4,5\\
    $q_2$ & $(\frac{1}{4},\frac{5\sqrt{7}}{28},\sqrt{\frac{5}{7}})$    & 4,5,6,7
  \end{tabular}
  \caption{Coordinates for vertices of $T$ in $\partial_\infty C$.}\label{tab:tvertices}
  \end{table}
The claims in the last column of the table follow from the results in
section~\ref{sec:combinatorics}, but they can also be checked directly
from their $t$-coordinates and the explicit expressions for $f_j$, see
Table~\ref{tab:equations2}.

From the equations for $\hB_2$, $\hB_4$ and $\hB_7$, one deduces
explicit parametri\-za\-tions for the three sides of $T$. For $\hB_4$
(and $\hB_5$), we get
  \begin{equation}\label{eq:b4c}
    (\frac{1}{16}(9-7t_3^2)\ ,\ \frac{1}{16}\sqrt{175-130t_3^2-49t_3^4\ },\ t_3),
  \end{equation}
for $t_3$ between $-\sqrt{5/7}$ and $\sqrt{5/7}$. This gives a
parametrization for $\tau_0$.

Here and in what follows, we write $f_j$, $j=1,\dots,8$ for the
equation for $\hB_j\cap C_\infty$ given in Table~\ref{tab:equations2}.
The parametrization for $\hB_2$ can be obtained by writing out the
resultant of $f_2$ and $t_1^2+t_2^2+t_3^2-1$ with respect to $t_1$,
which has degree 2 in $t_2$.  Using the quadratic formula, we get
  $$
  t_2 = \phi(t_3) = \frac{b(t_3) - (245 t_3 + 87\sqrt{5}\sqrt{7})\sqrt{d(t)}/35}{a(t_3)},
  $$
  where
  \begin{eqnarray*}
    & a(t) =  968 + 136\sqrt{5}\sqrt{7}t_3 + 320 t_3^2  \\  
    & b(t) = - 39\sqrt{5}t_3^3 + 80\sqrt{7}t_3^2 + 108\sqrt{5}t_3 - 55\sqrt{7}  \\
    & d(t) = - ( 1715 t_3^4 + 385 \sqrt{5}\sqrt{7}t_3^3 + 1750 t_3^2 + 175\sqrt{5}\sqrt{7}t_3 ).
  \end{eqnarray*}
  One then takes
  $$
  t_1 = \sqrt{1 - \phi(t_3)^2 - t_3^2},
  $$ 
and one checks that this parametrization is well-defined for $t_3$ in
the interval $[-\sqrt{5/7},0]$, which corresponds to the arc between
$q_1$ and $p_2$ of the triangle $T$. This gives a parametrization for
$\tau_2$.

We give the above explicit formulas mainly because there are two
solutions to the quadratic equation, so we need to select one.  The
parametrization for $\hB_7$ is obtained from the one for $\hB_2$
simply by changing $t_3$ into $-t_3$.
The latter property and the fact that the two paths on $\hB_2$ and
$\hB_7$ are parametrized by $t_3$ implies that these arc only
intersect along $t_3=0$, which corresponds to their common endpoint
$p_2$.

In order to prove that $\tau$ is embedded, it is enough to check that
the image of $\tau_0$ and $\tau_2$ intersect only in $q_1$ (the
corresponding property for $\tau_0$ and $\tau_1$ follows by symmetry).
The quickest way to show this is to compute a Groebner basis
for the ideal generated by $f_2$, $f_4$ and
$g(t_1,t_2)=t_1^2+t_2^2+t_3^2-1$, and to check that the corresponding
system has a unique solution, corresponding to $q_1$, or in other
words $(t_1,t_2,t_3)=(1/4,5\sqrt{7}/28,-\sqrt{5/7}).$

\begin{rmk}\label{rmk:quadrant}
  The path $\tau$ bounds two disks in $\partial_\infty C\simeq S^2$,
  only one of which is contained in the first quadrant $t_1,t_2>0$
  (this is the triangle $T$ that appears in
  section~\ref{sec:combinatorics}).
\end{rmk}

\subsection{Proof of Proposition~\ref{prop:triangle}} \label{sec:proofTriangle}

We denote by $T$ the (closure of) the component of the complement of
$\tau$ in $\partial_\infty C$ that is contained in the quadrant
$t_1,t_2>0$ in the coordinates of section~\ref{sec:proofdeftau4} (see
Remark~\ref{rmk:quadrant}). It is easy to see that the other component of its complement
is not contained in $U$, the difficult part is to show:
\begin{prop}\label{prop:properly}
  $T$ is properly contained in $U$.
\end{prop}

\begin{pf}
We first check that points on the boundary of $T$ are precisely on the
bisectors we think they are on (according to the incidence pattern
already mentioned in section~\ref{sec:combinatorics}). This can be
done by finding intersection points of pairs of curves corresponding
to the intersection of $\partial_\infty C$ with $\hB_j$, $\hB_k$,
$j\neq k$, which amounts to solving a system of equations, for
instance by using Groebner bases.

As an example, $\hB_1\cap \hB_2\cap \partial_\infty C$ has precisely
two points. One is $q_1$, and the other one is given approximately by
  $$(0.88541680, 0.03241871, -0.46366596).$$ It is easy to check that
this point is not in $T$.

With such verifications, one checks that the $\hB_j$ intersect $T$
only on its boundary, and only in a predicted fashion: the vertices
are on four bisectors, points in $]p_2,q_1[$ lie only in $\hB_2$ and
no other $\hB_k$, points in $]p_2,q_2[$ lie only in $\hB_7$ and in
no other $\hB_k$, points in $]p_1,p_2[$ lie on on $\hB_4$ and
$\hB_5$ and no other $\hB_k$.

We now rule out the possibility that some $\hB_j$ may have a
connected component contained in the interior of $T$. If that were the
case, then (the restriction to $\partial_\infty C$ of) $f_j$ would
have a critical point in the interior of $T$.

\noindent
\textsc{Claim:} no $f_j$ has a critical point in the interior of $T$.

A definite list of the critical points of $f_j$ can be obtained by using
Lagrange multipliers; the critical points for $f_j$ are the solutions
of the system
  \begin{equation}\label{eq:lagrange}
  \left\{
  \begin{array}{l}
  \nabla f_j = \lambda \nabla g\\
  g=0
  \end{array}\right.,
  \end{equation}
where $g(t)=t_1^2+t_2^2+t_3^2-1$. We only treat an example
representative of the difficulties, namely $f_2$. In that case, the
system~\eqref{eq:lagrange} reads
  $$
  \left\{\begin{array}{l}
  27/40 +  3\sqrt{7}t_3/8\sqrt{5} + (9/5-2\lambda)t_1 = 0 \\
  - 33 t_3/8\sqrt{5} - (3/8)\sqrt{7} + (9/5-2\lambda)t_2  = 0\\
  - 33 t_2/8\sqrt{5} + 3\sqrt{7}t_1/8\sqrt{5} - 3\sqrt{7}/2\sqrt{5}-2\lambda t_3 = 0\\
  t_1^2+t_2^2+t_3^2=1
  \end{array}\right.
  $$ 
This system can easily be solved using standard Groebner basis
techniques.

It has precisely four real solutions, for which $t_3$ is equal to $0$,
$-\sqrt{5/7}$, or one of the two real roots of the polynomial
  $$
  140t_3^4 + 28\sqrt{5}\sqrt{7}t_3^3 - 49t_3^2 - 20\sqrt{5}\sqrt{7}t_3 - 35, 
  $$ 
which are given approximately by $-0.50306965$ and $0.84223313$.  Only
one of the corresponding critical points lies in the first quadrant
$t_1,t_2>0$, and it corresponds precisely to $q_1$, which is not in
the interior of $T$.

For concreteness, we draw two projections of the $2$-sphere
$\partial_\infty C$, the triangle $T$ and the critical points of the
equations on Figure~\ref{fig:criticalpoints}. No critical points lies
in the interior of $T$, and the only critical points on the boundary
are $q_1$ (which is critical for $f_3$) and $q_2$ which is critical
for $f_6$). A couple of critical points may appear dubious on the
picture. One of them is $(x,y,t)=(1,0,0)$, which is critical for $f_1$
and $f_8$; one can easily check that it is not in $E$ by checking a
few inequalities (it is in fact only close to $p_2$, which has
approximate coordinates $(0.98863636,0.15032678,0)$, see
Table~\ref{tab:tvertices}).

Another pair of critical points are dubious only in
$(x,t)$-projection; on part~(b) of Figure~\ref{fig:criticalpoints},
they clearly appear outside the triangular region corresponding to
$T$.
  \begin{figure}
    \centering
    \subfigure[$(x,y)$-projection]{\epsfig{figure=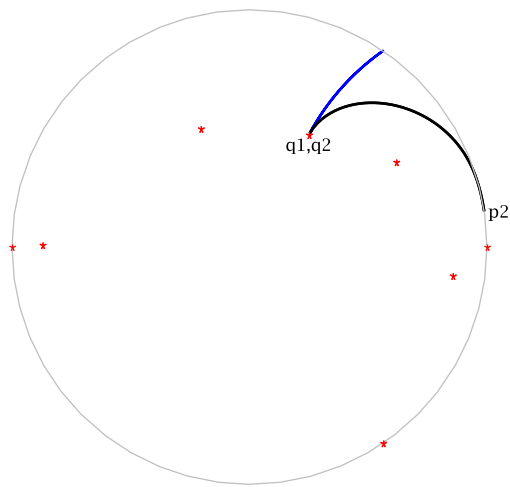, width=0.55\textwidth}}
    \subfigure[$(x,t)$-projection]{\epsfig{figure=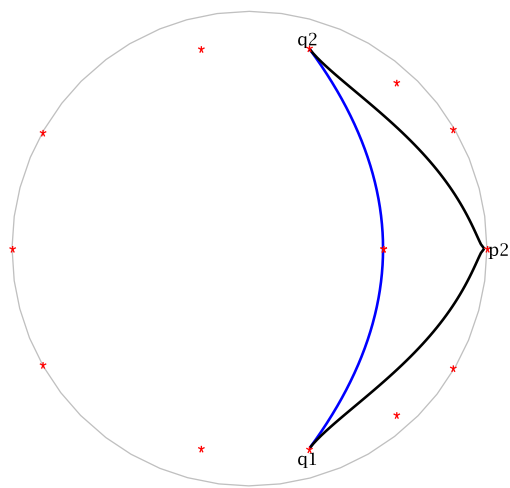, width=0.55\textwidth}}
    \caption{The critical points of the equations are outside $T$, in
      projection onto two coordinate
      planes.}\label{fig:criticalpoints}
  \end{figure}
\end{pf}

\begin{prop}
  The intersection $T\cap G_2^2 T$ is empty. 
\end{prop}

\begin{pf}
  We show a stronger statement, namely we show that $\partial_\infty
  C\cap G_2^2 \partial_\infty C$ consists of precisely two points that
  are not in $T$. We use the same coordinates for $C$ (and
  $\partial_\infty C$) as above, write $G_2^2G_1^{-1} p_0$ and
  $G_2^2G_2^{-1}G_3 p_0$ in terms of the basis given by the columns
  of~\eqref{eq:basechange2}, and write the equation of the
  intersection of $C$ with $G_2^2 C$, which is simply
  $$
  -12t_2t_3/\sqrt{5}.
  $$ 
  This gives $(1,0,0)$ and $(-1,0,0)$ as the only intersection
  points on $\partial_\infty$. None of these two points is in the
  Dirichlet domain, a fortiori they are not in $T$.
\end{pf}

\bibliographystyle{plain} \bibliography{biblio}

\begin{flushleft}
  \textsc{Martin Deraux\\
  Institut Fourier, Universit\'e de Grenoble 1, BP 74, Saint Martin d'H\`eres Cedex, France}\\
  \verb|deraux@ujf-grenoble.fr|
\end{flushleft}
\begin{flushleft}
  \textsc{Elisha Falbel\\
  Institut de Math\'ematiques,
  Universit\'e Pierre et Marie Curie,
  4 place Jussieu,
  F-75252 Paris, France}\\
  \verb|falbel@math.jussieu.fr|
\end{flushleft}

\end{document}